\documentclass[11pt]{article}

\usepackage{amsmath}
\usepackage{amssymb}
\usepackage{latexsym}
\usepackage{amsmath, amsfonts,amssymb, amsthm, euscript,makeidx,color,mathrsfs}
\usepackage{delarray}
\usepackage{tabularx}
\usepackage{amsbsy}
\usepackage{graphicx}
 \usepackage{layout}
\usepackage{booktabs}
\usepackage{cite}

\newtheorem{definition}{Definition}
\newtheorem{theorem}{Theorem}
\newtheorem{corollary}{Corollary}
\newtheorem{lemma}{Lemma}
\newtheorem{remark}{Remark}

\newtheorem{assumption}{Assumption}
\def\qed{ \ \vrule width.2cm height.2cm depth0cm\smallskip}

\newcommand{\hP}{\hat\dbP}

\newcommand{\ba}{\begin{array}}
\newcommand{\ea}{\end{array}}
\newcommand{\be}{\begin{equation}}
\newcommand{\ee}{\end{equation}}
\newcommand{\bea}{\begin{eqnarray}}
\newcommand{\eea}{\end{eqnarray}}
\newcommand{\beaa}{\begin{eqnarray*}}
\newcommand{\eeaa}{\end{eqnarray*}}

\def\a{\alpha}
\def\b{\beta}

\def\d{\delta}
\def\e{\varepsilon}

\def\l{\lambda}
\def\m{\mu}

\def\si{\sigma}
\def\t{\tau}

\def\th{\theta}
\def\o{\omega}

\def\D{\Delta}

\def\L{\Lambda}

\def\Om{\Omega}
\def\cF{{\cal F}}
\def\cG{{\cal G}}

\def\cQ{{\cal Q}}

\def\hC{\mathbb{C}}

\def\hE{\mathbb{E}}
\def\hF{\mathbb{F}}

\def\hP{\mathbb{P}}
\def\hQ{\mathbb{Q}}
\def\hR{\mathbb{R}}

\def\hX{\mathbb{X}}

\def\sB{\mathscr{B}}

\def\sP{\mathscr{P}}

\def\sU{\mathscr{U}}

\def\no{\noindent}

\def\ss{\smallskip}
\def\ms{\medskip}
\def\bs{\bigskip}
\def\q{\quad}
\def\qq{\qquad}

\def\cd{\cdot}

\def\lan{{\langle}}
\def\ran{{\rangle}}

\def\be{\begin{equation}}
\def\bel{\begin{equation}\label}
\def\ee{\end{equation}}
\def\bt{\begin{theorem}}
\def\bcd{\begin{condition}}
\def\ecd{\end{condition}}
\def\et{\end{theorem}}
\def\bc{\begin{corollary}}
\def\ec{\end{corollary}}
\def\bde{\begin{definition}}
\def\ede{\end{definition}}
\def\bl{\begin{lemma}}
\def\el{\end{lemma}}
\def\bp{\begin{proposition}}
\def\ep{\end{proposition}}
\def\br{\begin{remark}}
\def\er{\end{remark}}
\def\ba{\begin{array}}
\def\ea{\end{array}}
\def\ed{\end{document}}

\def\square#1{\vbox{\hrule\hbox{\vrule height#1%
     \kern#1\vrule}\hrule}}
\def\rectangle#1#2{\vbox{\hrule\hbox{\vrule height#1%
     \kern#2\vrule}\hrule}}
\def\qed{\hfill \vrule height7pt width3pt depth0pt}

\def\sqr#1#2{{\vcenter{\vbox{\hrule height.#2pt
              \hbox{\vrule width.#2pt height#1pt \kern#1pt \vrule width.#2pt}
              \hrule height.#2pt}}}}

\def\qed{ \hfill \vrule width.25cm height.25cm depth0cm\smallskip}
\newcommand{\dfnn}{\stackrel{\triangle}{=}}
\newcommand{\basa}{\begin{assumption}}
\newcommand{\easa}{\end{assumption}}

\newcommand{\bas}{\begin{assum}}
\newcommand{\eas}{\end{assum}}

\def\liminf{\mathop{\underline{\rm lim}}}

\def\lan{\mathop{\langle}}
\def\ran{\mathop{\rangle}}

 \def\cd{\cdot}

\def\as{\hbox{\rm -a.s.{ }}}

\def\dis{\displaystyle}

\def\1{{\bf 1}}

\def\:{\!:\!}

at 9pt
\begin{document}
\newtheorem{thm}{Theorem}[section]
\newtheorem{lem}[thm]{Lemma}
\newtheorem{cor}[thm]{Corollary}
\newtheorem{prop}[thm]{Proposition}
\newtheorem{rem}[thm]{Remark}
\newtheorem{eg}[thm]{Example}
\newtheorem{defn}[thm]{Definition}
\newtheorem{assum}[thm]{Assumption}

\renewcommand {\theequation}{\arabic{section}.\arabic{equation}}
\def\thesection{\arabic{section}}

\title{\bf  Kyle-Back Equilibrium Models and  \\
Linear Conditional Mean-field SDEs}
\author{
Jin Ma\footnote{Department of Mathematics, University of Southern Califorlia, Los Angeles, 90089, USA. Emial: jinma@usc.edu. This author is supported in part by US NSF grant DMS-1106853.}, ~{Rentao Sun}\thanks{\noindent
Department of Mathematics, University of Southern California, Los
Angeles, CA 90089, USA. E-mail: rentaosu@usc.edu.},
~and~ Yonghui Zhou\footnote{(Corresponding author) School of  Big Data and Computer Sciences, Guizhou Normal University, Guiyang, 550001, P.R. China. Email: yonghuizhou@gznu.edu.cn.  This author is supported by Chinese Scholarship Council \#20138525118 and  Chinese NSF grant \#11161011.} }

\date{}
\maketitle
\vspace{-3mm}
\begin{abstract}
In this paper we study the Kyle-Back strategic insider trading equilibrium model in which the insider has an instantaneous information on an asset, assumed to follow an Ornstein-Uhlenback-type dynamics that allows possible influence by the market price. Such a model exhibits some further interplay between insider's information and the market price, and it is the first time being put into a  rigorous mathematical framework of the recently developed {\it conditional mean-field}
stochastic differential equation (CMFSDEs). With the help of the ``reference probability measure" concept in filtering theory, we shall 
first prove a general well-posedness result for a class of linear CMFSDEs, which is new in the literature of both filtering theory and mean-field SDEs, and will be the foundation for the underlying strategic equilibrium model. Assuming  some further Gaussian 
structures of the model, 
we find  a closed form of optimal intensity of trading strategy as well as the dynamic pricing rules. We shall also substantiate the well-posedness of the resulting optimal closed-loop system, whence the existence of Kyle-Back equilibrium. Our result recovers many existing results as special cases.
\end{abstract}

\vfill \bs

\no{\bf Keywords.}  Strategic insider trading,  Kyle-Back equilibrium, conditional mean-field SDEs, reference measures, optimal closed-loop system. 

\bs

\no{\it 2000 AMS Mathematics subject classification:} 60H10, 91G80;
60G35, 93E11.



\bigskip

\section{Introduction}
In his seminal paper, A.S. Kyle \cite{K85}  first proposed a sequential equilibrium model of   asset pricing with asymmetric information. The model was then extended by K. Back \cite{B92} to the continuous time version, and has since known as the {\it Kyle-Back
strategic insider trading equilibrium model}. Roughly speaking, in such a model it is assumed that there are
two types of traders in a risk neutral market: one informed (insider) trader vs. many uninformed (noise) traders. The insider ``sees"
both (possibly future)  value of the fundamental asset as well as its market value, priced by the market makers, and acts
strategically in a non-competitive manner. The noise traders, on the other hand,
act independently with only market information of the asset. Finally, 
the market makers set the price of the asset, in a Bertrand competition fashion, based on the historical information of  the collective market actions of all traders, without knowing 
identify the insider. The so-called {\it Kyle-Back equilibrium} is a closed-loop system in which the insider maximizes his/her expected return 
in a {\it market efficient} manner (i.e., following the given market pricing rule).

There has been a large number of literature on this topic. We refer to, e.g.,  \cite{B92, B98, CS, Cho, CLS,  Dani, FV, GM, HS, K85} and the references therein for both discrete and continuous time models. It is noted, however, that in most of these works only the case of static information is considered, that is, it is assumed that information that the insider could observe is ``time-invariant", often as the
fundamental price at a given future moment. Mathematically, this amounts to saying that the insider has the knowledge of a
given random variable whose value
%
cannot be detected from the market information at current time.
It is often assumed that the system has a certain Gaussian structure (e.g., the future price is a Gaussian random variable),
so that the optimal strategies can be calculated explicitly.
The situation will become more complicated when the fundamental price progresses as a stochastic process $\{v_t, t\ge 0\}$
and the insider is able to observe the price process dynamically in a ``non-anticipative" manner. The  asymmetric information nature of the problem has conceivably led to the use of {\it filtering} techniques 
in  the study of the Kyle-Back model, and we refer to, e.g., \cite{ABO12} and
\cite{BHMO12} for the static information case, and to, e.g., \cite{CS} and \cite{Dani} for the dynamic information case.
It is noted that
in \cite{BHMO12} it is further assumed that the actions of noise traders may have some ``memory", so that the observation process
in the filtering problem is driven by a fractional  Brownian motion, adding  technical difficulties in a different aspect. We also note
that the Kyle-Back model has been continuously extended in various directions. For example, in a static information setting, \cite{CD-F} recently
considered the case when noise trading volatility is a stochastic process, and in the dynamic information case \cite{CCD11, CCD13a, CCD13b} studied the Kyle-Back equilibrium for the defaultable underlying asset via dynamic Markov bridges, exhibiting
further theoretical potential  of the problem.

In this paper we are interested in a generalized Kyle-Back equilibrium model in a dynamic information setting, in which the asset dynamics is of the form of an Ornstein-Uhlenback SDE whose drift also  contains market sentiment (e.g.,
supply and demand, earning base, etc.), quantified by
 the market price. The problem is then naturally imbedded into a (linear) filtering problem in which both  state and observation dynamics contain the filtered signal process (see \S2 for details). We note that such a structure is not covered by the existing 
 filtering theory,  and thus it is interesting in its own right.
In fact, under the setting of this paper the signal-observation dynamics form a  ``coupled" (linear) conditional mean-field stochastic
differential equations (CMFSDEs, for short) whose well-posedness, to the best of our  knowledge, is new.

The main objective of this paper is thus two-fold. First, we shall look for a rigorous framework on which the well-posedness of the underlying CMFSDE can be established.
The main device of our approach is the ``reference probability measure" that is often seen in the nonlinear filtering
theory (see, e.g., \cite{Zeit}). Roughly speaking, we give the observable market movement a ``prior" probability distribution  so that it is a Brownian motion that is independent of the martingale representing the aggregated trading actions of the noisy traders, and we then prove
that the original signal-observation SDEs have a weak solution. More importantly, we shall prove that the uniqueness holds
among all weak solutions whose laws are absolutely continuous with respect to the reference probability measure. We should 
note that such
 a uniqueness in particular resolves a long-standing issue on the Kyle-Back equilibrium model: the identification of the total traded volume movements and the innovation process of the corresponding filtering problem, which has been argued either only heuristically
or by economic instinct in the literature (see, e.g.,  \cite{ABO12}).
 The second goal of this paper is to identify the Kyle-Back equilibrium, that is, the optimal closed-loop system for this new
 type of partially observable optimization problem.  
Assuming a Gaussian structure and the linearity of the CMFSDEs, 
we give the explicit solutions to the insider's trading intensity and justify the well-posedness of the closed-loop system (whence the ``existence" of the equilibrium).
These solutions in particular cover many existing results as special cases.

%


The rest of the paper is organized as follows. In section 2 we give the preliminaries of the Kyle-Back equilibrium model and formulate the strategic insider trading problem. In section 3 we formulate a general form of linear CMFSDE and introduce the notion
of its solutions and their uniqueness. We state the main well-posedness result, calculate its explicit solutions in the case of deterministic coefficients, and discuss an important extension to the unbounded coefficients case. Section 4 will be devoted to the 
proof of the main well-posedness theorem. In section 5 we characterize the
optimal trading strategy, and give a first order  necessary condition for the optimal intensity, and finally in section 6 we focus on
the synthesis analysis, and validate 
the Kyle-Back equilibrium. In some spacial cases, we give the closed-form solutions, and compare them to the existing results.
%

\section{Problem Formulation}
\setcounter{equation}{0}

In this section we describe a continuous time Kyle-Back equilibrium model that will be investigated in this paper,
as well as related technical settings.

We begin by assuming that all randomness of the market comes from
a common complete probability space $(\Om,\cF, \hP)$ on which is
defined 2-dimensional Brownian motion $B=(B^v, B^z)$, where $B^v=\{B^v_t:t\ge 0\}$
represents the noise of the fundamental value dynamics, and $B^z=\{B^z_t:t\ge 0\}$ represents
the collective action of the noise traders.
For notational clarity, we denote $\hF^v=\{\cF^{B^v}_t:t\ge 0\}$ and
$\hF^z\dfnn\{\cF^{B^z}_t:t\ge0\}$ to be the filtrations generated by
$B^v$ and $B^z$, respectively, and denote $\hF=\hF^v\vee\hF^z$,
with the usual $\hP$-augmentation such that it satisfies the {\it
usual hypotheses} (cf. e.g., \cite{prot}).

Further, throughout the paper we will denote, for a  generic Euclidean space
$\hX$, regardless of its dimension, $\lan\cd,\cd\ran$ and $|\cd|$ to be its inner product and norm, respectively.
We denote the space of continuous functions defined on $[0,T]$ with the usual sup-norm by $\hC([0,T];\hX)$, and we shall make use
of the following notations:
\begin{itemize}
\item{} For any sub-$\si$-field $\cG\subseteq\cF_T$ and $1\le p<\infty$,
$L^p(\cG;\hX)$ denotes the space of all $\hX$-valued, $\cG$-measurable
random variables $\xi$ such that $\hE|\xi|^p<\infty$. As usual, $\xi\in L^\infty
(\cG;\hX)$ means that it is $\cG$-measurable and bounded.

\item{} For $1\le p<\infty$, $L^p_\hF([0,T];\hX)$ denotes the space of all
$\hX$-valued, $\hF$-progressively measurable processes $\xi$ satisfying
$\hE\int_0^T|\xi_t|^pdt<\infty$. The meaning of $L^\infty _\hF([0,T];\hX)$
is defined similarly.
\end{itemize}
Throughout this paper we assume that all the processes are 1-dimensional, but higher dimensional cases can be
easily deduced without substantial difficulties. Therefore, we will often drop $\hX(=\hR)$ from the notation. Also, throughout the paper we shall denote all ``$L^p$-norms" by
$\|\cd \|_p$, regardless its being $L^p(\cG)$ or $L^p_{\hF}([0,T])$, when the context is clear.

Consider a given stock whose fundamental value (or its return) is $V=\{V_t:t\ge0\}$ traded on a finite time interval $[0,T]$. There are three types of agents in the market:
\begin{itemize}
\item[ (i)] {\it The insider}, who directly observes the realization of the value $V_t$ at any time $t\in [0,T]$ and submits his order $X_{t}$ at $t\in[0,T]$.

\item[(ii)] {\it The noise traders}, who have no direct information of the given asset, and (collectively) submit an order $Z_t$ at $t\in[0,T]$ in the form
\begin{equation}\label{z}
Z_t=\int_0^t\sigma^z_{t}dB^z_{t},\qq t\ge 0,
\end{equation}
where $\si^z=\{\sigma^z_t: t \ge 0\}$ is a given continuous deterministic function satisfying $\si^z_t>0$.

\item[(iii)] {\it The market makers},  who {observe} only the total traded volume
\begin{equation}
\label{y}
Y_t=X_{t}+Z_t,\qq t\ge 0,
\end{equation}
and set the market price of the underlying asset at each time $t$, denoted by $P_t$, based on the observed information 
${\cal F}_t^Y\dfnn \si\{Y_s, s\leq t\}$. 
We denote  $S_t\dfnn \hE[(V_t-P_t)^2]$ to be the measure of the discrepancy between  market price and fundamental value of the asset. 
\end{itemize}

We now give a more precise description of the {two} main ingredients in the model above:
the dynamics of market price $P$ and the fundamental value $V$. First,  in the same spirit of the original Kyle-Back model, we can assume that the market price $P$ is the result of a Bertrand-type competition among the market makers (cf. e.g., \cite{B98}), and therefore should be taken as $P_t={ \hE}[V_t|{\cal F}_t^Y]$, at each time $t\ge 0$. 
Mathematically speaking, this amounts to saying that  the market price is set to minimize the error $S_t$, among
all 
$\cF^Y_t$-measurable random variables in $L^2(\Om)$, hence the ``best estimator" that the market maker is able to choose given the information $\cF^Y_t$. It is easy to see that under such a choice one must have  $\hE[(V_t-P_t)Y_t]=0$, 
that is, the market makers should expect a zero profit at each time $t\in [0,T]$.
It is worth noting that, however, unlike the static case (i.e., $V_t\equiv v$), the process $P$ is no longer a $(\hP,\hF^Y)$-martingale in general.

Next, in this paper we shall also assume that the dynamics of the value of the stock $V=\{V_t\}$ takes the form of an It\^o process:
$dV_t=F_tdt+\si^v_tdB^v_t$, $t\ge 0$ (this would easily be the case if, e.g., the interest rate is non-zero).
Furthermore, we shall assume that the drift  $F_t=F(t,V_t, P_t)$, $t\ge0$. Here, the dependence { of $F$}
on the market price $P_t$ is based on the following rationale:  the value of the stock is often affected by factors such as supply
and demand, the earnings base (cash flow per share), or more generally, the market sentiment, which  all depend on the market
price of the stock.
Consequently, taking the Gaussian structure into consideration, in what follows we shall assume that the process $V$ satisfies the following
linear SDE:
  \bea
\label{V}
\left\{\ba{lll}
dV_t=(f_tV_t+g_tP_t+h_t)dt+\sigma^v_{t}d B^v_{t} =(f_tV_t+g_t{ \hE}[V_t|{\cal F}_t^Y]+h_t)dt+\sigma^v_{t}d B^v_{t}, \q t\ge 0,
\\
V_0\sim N(v_0,s_0).
\ea\right.
\eea
where the functions $f_t, g_t, h_t$ and $\si_t^v$ are all  deterministic continuous differentiable  functions  with respect to time $t\in [0, T]$,
and  $N(v_0, s_0)$ is a  normal random variable with mean $v_0$ and standard deviation $s_0$.
%
%

 Continuing, given the Gaussian structure of the dynamics, it is reasonable to
assume that the insider's optimal trading strategy (in terms of ``number of shares") is of the form (see, e.g., \cite{ABO12, B92, B98, Dani, K85}):
\bea\label{x0}
dX_t=\beta_t(V_t-P_t)dt, \qq t\ge 0,
\eea
 where $\b_t>0$ is a deterministic continuous differentiable  function  with respect to time $t$ in  [0, T),  known as the {\it insider trading intensity}.
 Consequently, it follows from (\ref{y}) that the total traded volume process observed by the  market makers can be expressed as
 \bea
 \label{y1}
 dY_t=\beta_t(V_t-P_t)dt+\sigma_t^zdB_t^z=\beta_t(V_t-\hE[V_t|\cF^Y_t])dt+\si^z_tdB^z_t, \qq t\ge 0.
 \eea

 We note that SDEs (\ref{V}) and (\ref{y1}) form a (linear) conditional mean-field SDE (CMFSDE),
 which is beyond the scope of the  traditional filtering theory. Such SDE
have been studied in \cite{BLM} and \cite{CZ} in  general nonlinear forms, but none of them covers the one in this  
form.
In fact, if we further allow the function $h$ in (\ref{V}) to be an $\hF^Y$-adapted process, as many stochastic control problems do, then (\ref{V}) and (\ref{y1}) would become a fully convoluted CMFSDE whose well-posedness, to the 
best of our knowledge,  has not been studied in the literature, even in the linear form.  We should mention that the equation (\ref{y1}) in the case when $V_t\equiv v$ was already noted in \cite{ABO12} and  \cite{BHMO12}, but without
addressing the uniqueness of the solution. In the next sections we shall establish a mathematical framework so these SDEs can be studied rigorously.  


Given the dynamics (\ref{V}) and (\ref{y1}), our main purpose now is to find an optimal trading intensity $\beta$  for the insider to maximize his/her expected  wealth, whence the Kyle-Back equilibrium. More specifically, denote the wealth process of the insider by $W=\{W_t:t\ge 0\}$, and assume that the
strategy is self-financing (cf. e.g., \cite{B09}), then
 the total wealth of the insider over time duration $[0,T]$, based on the market price made by the market makers, should
be
\bea
\label{w}
W_{T}=\int_0^T X_{t}dP_t=X_TP_T-\int_0^T\beta_{t}(V_t-P_t)P_tdt=\int_0^T\beta_{t}(V_t-P_t)(P_T-P_t)dt.
\eea
Here in the above we used a simple integration by parts and definition (\ref{x0}). Thus the optimization problems can
be described as
\bea
\label{j0}
\sup_{\b}\hE[W_T]\dfnn \sup_{\b}J(\beta)=\sup_{\b}\int_0^T\beta_{t}\hE[(V_t-P_t)(P_T-P_t)]dt.
\eea
%
%

%

 \begin{rem}
\label{remark02}{\rm
We should remark that the simple form of optimization problem (\ref{j0}) is due largely to the linearity of the dynamics (\ref{V}) and (\ref{y1}),
as well as the Gaussian assumption on the initial state $v$. These lead to a Gaussian structure, whence the trading strategy (\ref{x0}).
The general
nonlinear and/or non-Gaussian Kyle-Back model requires further study of CMFSDE and associated filtering problem,
and one should seek optimal control from a larger class of ``admissible controls". In that case the first order condition studied in
this paper will become a Pontryagin type stochastic maximum principle
(see, for example, \cite{BLM}), and the solution is expected to be much more involved. We
will address such general problems in our future publications.
\qed}
\end{rem}

We end this section by noting that the main idea for solving the CMFSDE is to introduce the  so-called
{\it reference probability space} in nonlinear filtering literature (see, e.g.,  \cite{Zeit}), which can be  described as follows.
\begin{assum}
\label{assump1}
%
There exists  {a  probability space $(\Om^0, \cF^0, \hQ^0)$}  on which the process $(B^v_t, {Y_t})$, $t\in [0,T]$,
is a 2-dimensional continuous martingale, where $B^v$ is a standard Brownian motion and $Y$ is the observation process  with quadratic variation $\lan Y\ran_t=\int_0^t(\sigma_s^z)^2ds$. The probability measure $\hQ^0$ will be referred to
as the {\it reference measure}.
\end{assum}
We remark that Assumption \ref{assump1} amounts to saying that we are giving a {\it prior} distribution to the price process
$Y=\{Y_t:t\ge0\}$ that
the market maker is observing, which is not unusual in statistical modeling, and will facilitate the discussion greatly.
A natural example is the {\it canonical space}: $\Om^0\dfnn C_0([0,T];\hR^2)$, the
 space of all 2-dimensional continuous functions null at zero;
$\cF^0\dfnn\sB(\Om^0)$; $\cF^0_t\dfnn \sB_t(\Om^0)=\si\{\o(\cd\wedge t): \o \in\Om^0\}$, $t\in[0,T]$; and $(B^v, Y)$ is the
canonical process. In the case {$\si^z\equiv 1$}, $\hQ^0$ is the Wiener measure.

\section{The Linear Conditional Mean-field SDEs}
\setcounter{equation}{0}

In this section we study the linear conditional mean-field SDEs (CMFSDE) (\ref{V}) and (\ref{y1}) that play an important role in this
paper. In fact, let us consider a slightly more general case that is useful in applications:
\bea
\label{sde4.1}
\left\{
\ba{lll}
dX_t=\{f_tX_t+g_t \hE [X_t|{\cal F}_t^Y]+h_t \}dt+ \si^1_td B^1_{t},\q X_0=v;\\
dY_t=\{H_t X_t + G_t \hE [X_t|{\cal F}_t^Y]\}dt+ \si^2_tdB_t^2,\q  Y_0=0,
\ea
\right.
\eea
where $B\dfnn (B^1, B^2)$ is a standard Brownian motion defined on a given probability space $(\Om, \cF, \hP)$, and
$v \sim N(v_0,s_0)$ is independent of $B$.
In light of Assumption \ref{assump1}, throughout this section we shall assume the following:
\begin{assum}
\label{assump2}
(i) The coefficients $f$, $g$, $\si^1$, $\si^2$, $G$, and $H$ are all deterministic, continuous functions; and $\si^i_t > 0$, $i=1,2$, for all $t\in[0,T]$;
%

(ii) there exists a probability space $(\Om^0, \cF^0, \hQ^0)$ on which the process $(B^1_t, Y_t)$, $t\in[0,T]$,
is a 2-dimensional continuous martingale, such that $B^1$ is a standard $\hQ^0$-Brownian motion, and
$\lan Y\ran_t=\int_0^t|\sigma_s^2|^2ds$, $t\in[0,T]$,  $\hQ^0$-a.s.;

(iii) the coefficient $h$ is an $\hF^Y$-adapted,  continuous process, such that
$$\hE^{\hQ^0}\Big[\sup_{0\le t\le T}|h_t|^2\Big]<\infty.
$$
\end{assum}

\begin{rem}
\label{rem1}
{\rm  The Assumption \ref{assump2}-(iii) amounts to saying that the process $h$ is defined on the reference
probability space $(\Om^0, \cF^0, \hQ^0)$, and adapted to the Brownian filtration $\hF^Y$, as we often see in the stochastic
control with partial observations (cf. \cite{Bens}).
\qed}
\end{rem}

\subsection{The General Result}

To simplify notations in what follows we shall assume that $\sigma^1=\sigma^2\equiv 1$. We first introduce two definitions of the solution
to  CMFSDE (\ref{sde4.1}).  Let $\sP(\hR)$ denote all probability measures on $(\hR, \sB(\hR))$, where
$\sB(\hR)$ is the Borel $\si$-field of $\hR$, and $\m\sim N(v_0,s_0)\in\sP(\hR)$ denotes the normal distribution with mean $v_0$
and variance $s_0$.
\begin{defn}
\label{weak sol}
Let $\m\in\sP(\hR)$ be given. An eight-tuple $(\Om,\cF,\hF,\hP; X, Y, B^1, B^2)$ is called a weak solution to CMFSDE (\ref{sde4.1})
with initial distribution $\m$ if

(i) $(B^1, B^2)$ is an $\hF$-Brownian motion under $\hP$;

(ii) $(X, Y, B^1, B^2)$ satisfies (\ref{sde4.1}), $\hP$-a.s.;

(iii) $X_0\sim \m$; and is independent of $(B^1, B^2)$ under $\hP$.
\qed
\end{defn}

\begin{defn}
\label{Q0weak}
A weak solution $(\Om,\cF,\hF,\hP; X, Y, B^1, B^2)$ is called a $\hQ^0$-weak solution {to  CMFSDE (\ref{sde4.1})} if

(i) there exists a probability measure $\hP^0$ on $(\Om^0, \cF^0)$, and
processes $(X^0, Y^0, B^{1,0}, B^{2,0})$ defined on $(\Om^0, \cF^0,\hP^0)$, whose law under $\hP^0$ is the same as
that of $(X, Y, B^1, B^2)$ under $\hP$; and

(ii) $\hP^0\sim \hQ^0$.
\qed
\end{defn}

In what follows for any $\hQ^0$-weak solution, we shall consider only its copy on the reference measurable space $(\Om^0, \cF^0)$, and we shall still denote the solution by $(X, Y, B^1, B^2)$.

The uniqueness of the solutions to  CMFSDE (\ref{sde4.1}) is a more delicate issue. In fact,
even the weak uniqueness (in the usual sense) for  CMFSDE (\ref{V}) and (\ref{y1}) is not clear. However,  we have a much better hope,
 at least in
 the linear case, for $\hQ^0$-solutions. We first introduce the following ``{\it $\hQ^0$-pathwise uniqueness}".
\begin{defn}
\label{unique}
The CMFSDE (\ref{sde4.1}) is said to have ``{\it $\hQ^0$-pathwise  uniqueness}" if for any two $\hQ^0$-weak solutions
$(\Om^0,\cF^0,\hF^0,\hP^i; X^i, Y^i, B^{1, i}, B^{2,i})$, $i=1,2$, such that

(i) $X^1_0=X^2_0$; and

(ii)  $\hQ^0\{(B^{1,1}_t, Y^1_t)= (B^{1,2}_t, Y^2_t), \forall t\in[0,T]\}=1$,

\no then it holds that $\hQ^0\{(X^{1}_t, B^{2,1}_t)= (X^2_t, B^{2,2}_t), \forall t\in[0,T]\}=1$, and $ \hP^1=\hP^2$.
\qed
\end{defn}

\begin{thm}
\label{well-posed}
{ Assume  that Assumption  \ref{assump2}} is in force, and  further that $h$ is bounded. Let $\m\sim N(v_0,s_0)$ be given. Then  CMFSDE (\ref{sde4.1}) possesses a weak solution
with initial distribution $\m$, denoted by $(\Om,\cF,\hF,\hP; X, Y, B^1, B^2)$.

Moreover, if we denote $P_t=\hE^\hP[X_t|\cF^Y_t]$, $t\in[0, T]$, then $P$ satisfies the following SDE:
\bea
\label{sde4.19}
\left\{\ba{lll}
dP_t=[(f_t +g_t )P_t +h_t ]dt+S_t H_t \{dY_t -[H_t+G_t] P_tdt\}, \q t\in[ 0,T], \ms\\\
P_0 =v_0.
\ea\right.
\eea
where $S_t=$Var$(P_t)$ satisfies the Riccati equation:
\bea
\label{SDE4.21}
dS_t=[1+2f_t S_t-H_t^2 S_t^2]dt, \q S_0= s_0.
\eea

Furthermore, the weak solution can be chosen as $\hQ^0$-weak solution, and the $\hQ^0$-pathwise uniqueness
holds.
\qed
\end{thm}

We remark that Theorem \ref{well-posed} does not imply that {CMFSDE (\ref{sde4.1})} has a strong solution, as not
every weak solution is a $\hQ^0$-weak solution. 
The proof of Theorem \ref{well-posed} is a bit lengthy, we shall defer it to next section. We nevertheless present a lemma below, which will be frequently used in our discussion,  so as to facilitate the
argument in the next section.

To begin with, we consider any filtered probability space  $(\Om,\cF, \hF, \hP)$ on which is defined a standard Brownian Motion
$(B^1_t, B^2_t)$. We assume that $\hF=\hF^{(B^1, B^2)}$. For any $\eta\in L^2_{\hF}([0,T])$ we define $L^\eta$ to be
the solution to the following SDE,
\be
\label{sde4.5}
dL_t=L_t \eta_tdB^2_t, \qq t\ge 0, \qq L_0=1.
\ee
In other words, $L^\eta$ is a local martingale in the form of the Dol\'eans-Dade stochastic exponential:
\bea
\label{Leta}
L^\eta _t=\exp\Big\{\int_0^t \eta_sdB^2_s-\frac12\int_0^t|\eta_s|^2ds\Big\}.
\eea
Next let  $\a\in L^2_\hF([0,T])$ and consider the following SDE:
\be
\label{Y}
dY_t= (\a_t+h(Y)_t)dt + dB^2_t, \qq Y_0=0,
\ee
where $h: [0,T]\times \hC([0,T])\mapsto \hR$ is ``progressively measurable" in the sense that, it is a measurable function such
that for each $t\in[0,T]$,
$h(y)_t =h(y_{\cd\wedge t})_t$ for $y\in \hC([0,T])$. (A simple case would be $h(y)_t=\tilde h(y_t)$, where $\tilde h$ is a measurable function.)
We should note that in general the well-posedness of SDE (\ref{Y})
is non-trivial without any specific conditions on $h$, but in what follows we shall assume {\it a priori} that (\ref{Y}) has a (weak) solution on some probability space $(\Om, \cF, \hP)$. We say that $h\in L^2_{\hF^Y}([0,T])$ if $h_t=h(Y)_t$, $t\in[0,T]$, such that $\hE\int_0^T|h(Y)_t|^2dt<\infty$. We have the following lemma.
\begin{lem}
\label{lem4.1}
Suppose that the SDE (\ref{Y}) has a solution $Y_t$, $t\in[0,T]$, for given $\a_t\in L^2_{\hF^{B^1}}([0,T])$ and $h\in L^2_{\hF^Y}([0,T])$ on some probability space $(\Om, \cF, \hP)$. Let $\b$  be given by
\be
\label{beta}
d\b_t= \a_t dt+ dB_t^2,\q t\ge 0, \q \b_0=0.
\ee
Assume further that $L^{-(\a+h)}$,
the solution to (\ref{sde4.5}) with $\eta=-(\a+h)$, is an $(\hF, \hP)$-martingale.
Then, for any $t \in [0, T]$, it holds that
\be
\label{sde4.17}
\hE^{\hP} [\a_t|\cF_t^Y]=\hE^{\hP} [\a_t|\cF_t^{\b}], \qq \forall t\in [0, T], \q \hP\as
\ee
\end{lem}

{\it Proof.}
Clearly, it suffices to prove
$\hE^\hP [\a_T|\cF_T^Y]=\hE^\hP [\a_T|\cF_T^{\b}]$, as the cases for $t<T$ are analogous. To this end, we define a new probability measure $\hQ$ on $(\Om, \cF_T)$ by
\beaa
\frac{d\hQ}{d\hP}\Big|_{\cF_T}=L^{-(\a+h)}_T,
\eeaa
where $L^{-(\a+h)}$ is the solution to the SDE (\ref{sde4.5}) with $\eta=-(\a+h)$, and it is a true martingale on $[0,T]$ by assumption.
By Girsanov Theorem,
 the process $(B^1, Y)$ is a standard Brownian motion on $[0,T]$ under $\hQ$.

Now define $\bar L_t =1/ L^{-(\a+h)}_t$, then $\bar L$ satisfies the following SDE on $(\Om, \cF_T, \hQ)$:
\be
\label{sde4.21}
d\bar L_t = \bar L_t(\a_t+h_t)dY_t, \qq t \in[0,T], \q \bar L_0=1.
\ee
Furthermore, by the Kallianpur-Striebel formula, we have
\be
\label{KS}
\hE^\hP [\a_T| \cF_T^Y]= \frac{\hE^\hQ[\a_T\bar L_T | \cF_T^Y]}{\hE^\hQ[\bar L_T | \cF_T^Y]}.
\ee
On the other hand, $\bar L$ has the explicit form:
\bea
\label{SDE}
\bar L_T&=&\exp\Big\{ \int_0^T [\a_t+h_t]dY_t - \frac{1}{2} \int_0^T [\a_t+h_t]^2dt\Big\}\\
&=&\exp\Big\{ \int_0^T h_tdY_t -\frac{1}{2}\int_0^T|h_t|^2dt+\int_0^T \a_tdY_t - \frac{1}{2} \int_0^T[ |\a_t|^2+2 h_t \a_t ]dt\Big\}
\nonumber\\
&\dfnn& \bar L^0_T\L_T,  \nonumber
\eea
where
\beaa
\bar L^0_T\dfnn \exp\Big\{\int_0^T h_tdY_t - \frac{1}{2} \int_0^T|h_t|^2dt \Big\}; ~\L_T\dfnn\exp\Big\{\int_0^T \a_tdY_t -
\frac{1}{2} \int_0^T[ [\a_t]^2+2 h_t \a_t ]dt\Big\}.
\eeaa
Note that $h$ is $\hF^Y$-adapted, so is $\bar L^0_T$. We derive from (\ref{KS}) that
\be
\hE^\hP [\a_T| \cF_T^Y]= \frac{\hE^\hQ[\a_T\L_T | \cF_T^Y]}{\hE^\hQ[\L_T | \cF_T^Y]}.
\ee

Now define  $Y^1_t=\int_0^t h_sds$. Since  $h$ is  $\hF^Y$-adapted, so is $Y^1$, and consequently
$\b_t=Y_t-Y^1_t$,
$t\ge0$ is $\hF^Y$-adapted. Moreover, since $(B^1,Y)$ is a standard Brownian motion under $\hQ$, and $\a$ is $\hF^{B^1}$-adapted, we conclude that $\a_t$ is independent of $Y_t$ under $\hQ$. Therefore, using integration by parts we obtain that
\bea
\label{Lambda}
\L_T&=& \exp\Big\{\int_0^T \a_td\b_t - \frac{1}{2} \int_0^T|\a_t|^2dt \Big\}\\
&=&\exp\Big\{ \a_T \b_T-\int_0^T \b_t d \a_t - \frac{1}{2}\int_{0}^{T}|\a_t|^2 dt\Big\}.
\nonumber
\eea
Since $\a$ is independent of $Y$  under $\hQ$, and $\b_t$, $t\in [0,T]$ is $ \cF_T^Y$-measurable, a Monotone Class argument shows that
$\hE^\hQ[\L_T | \cF_T^Y]$ is $\cF_T^{\b}$ measurable; and similarly, $\hE^\hQ[\a_T\L_T | \cF_T^Y]$ is also $\cF_T^{\b}$ measurable. Consequently $\hE^\hP [\a_T| \cF_T^Y]$ is $\cF_T^{\b}$ measurable, thanks to (\ref{KS}).

Finally, noting $\hF^{\b}\subseteq \hF^Y$ we have
\be
\hE^\hP [\a_T| \cF_T^Y]=\hE^\hP \{\hE^\hP [\a_T| \cF_T^Y]|\cF_T^{\b}\}=\hE^\hP [\a_T| \cF_T^{\b}],
\ee
proving the lemma.
\qed

\subsection{Deterministic Coefficient Cases}

An important special case is when {\it all} the coefficients in the linear CMFSDE (\ref{sde4.1}) are deterministic. In this case
we expect that the solution $(X, Y)$ is Gaussian, and it can be solved in a much more explicit way. The following  linear
CMFSDE will be useful in the study of insider trading equilibrium model in the latter half of the paper.
\bea
\label{v1}
\left\{\ba{lll}
dX_t=[f_tX_t+g_t\hE[X_t|{\cal F}_t^Y]+h_t]dt+\sigma^1_{t}d B^1_{t}, \q X_0=v;\ms\\
dY_t=H_t(X_t-\hE[X_t|{\cal F}_t^Y])dt+\sigma_t^2dB_t^2,\q  Y_0=0,
\ea\right.
\eea
where $v\sim N(v_0,s_0)$ and is independent of $(B^1, B^2)$ and all the coefficients are assumed to be deterministic.

\subsubsection{Bounded Coefficients Case}
In light of Theorem \ref{well-posed} let us introduce the following functions:
\bea
\label{kl}
k_t= \frac{H^2_tS_t}{|\si_t^2|^2},\q \q  \ l_t=\frac{H_tS_t}{\si_t^2}, \qq t\ge0.
\eea
where  $S$ is the solution to the following Riccati equation
\bea
\label{Riccati}
\frac{d S_t}{dt}=(\si_t^1)^2+2f_tS_t-l_t^2, \qq t\ge0, \qq S_0=s_0.
\eea

We have the following result.
\begin{prop}
\label{eu}
Let  Assumption \ref{assump2} be in force, and assume further that the process $h$ in (\ref{sde4.1}) is also a deterministic
and continuous function.
Let $(X,Y)$ be the solution of (\ref{v1}) on the probability space $(\Om, \cF, \hP)$, and denote $P_t=\hE^\hP[X_t|\cF^Y_t]$,
$t\ge0$. Then $X$ and $P$ have the following explicit {forms respectively}: for $t\ge0$, it holds $\hP$-a.s. that
\bea
\label{vmp}
X_t&=&P_t+\phi_1(t,0)\Big[v-v_0+\int_0^t\phi_1(0,r)(\sigma^1_{r}d B^1_{r}-l_rd Y_{r})\Big]; \\
\label{p1}
 P_t&=&\phi_2(t,0)\Big\{v_0+\int_0^t\phi_2(0,r)h_rdr+(v-v_0)\phi_3(t,0)\\
\nonumber &&+\int_0^t\sigma_r^1\phi_1(0,r)\phi_3(t,r)dB_r^1+
\int_0^t[\phi_2(0,r)l_r-\phi_1(0,r)\phi_3(t,r)l_r]dB_r^2\Big\},
\eea
where, for  $0\leq r\leq t$,
\bea
\label{phi123}
\left\{\ba{lll}
\dis \phi_1(t,r)=\exp\Big\{\int_r^t(f_u-k_u)du\Big\}; \q \phi_2(t,r)=\exp\Big\{\int_r^t(f_u+g_u)du\Big\};\ms\\
\dis \phi_3(t,r)=\int_r^t\phi_1(u,0)\phi_2(0,u)k_udu.
\ea\right.
\eea
 \end{prop}

 {\it Proof.} We first note that the SDE (\ref{v1}) is a special case of (\ref{sde4.1}) with $G=-H$. Then, following the same
 argument of Theorem \ref{well-posed} one can show that
when $\si^1>0$ and $\si^2>0$ are not equal to 1, the SDE (\ref{sde4.19}) for the process  $P_t=\hE^\hP[X_t|{\cal F}_t^Y]$
reads
\bea
\label{p2} dP_t=[(f_t+g_t)P_t+h_t]dt+\frac{H_tS_t}{(\si_t^2)^2}dY_t, \qq P_0=v_0,
\eea
and $S_t$ satisfies a Riccati equation
\bea
\label{s1}
\frac{dS_t}{dt}=(\si_t^2)^2+2f_tS_t-\Big[\frac{H_tS_t}{\si_t^2}\Big]^2,\qq  S_0=s_0.
\eea

Now applying the Girsanov transformation we can define a new probability measure  $ \hQ$ under which $(B^1, Y)$ is a continuous martingale, such that $B^1$ is a standard Brownian motion, and $d\lan Y\ran_t=|\si^2_t|^2dt$. Then, under $\hQ$, the
 dynamic of $V_t\dfnn X_t-P_t$ can be written as
\beaa
dV_t&=&\Big[f_t-\frac{H^2_tS_t}{(\si_t^2)^2}\Big]V_tdt +\sigma^1_{t}d B^1_{t}-\frac{H_tS_t}{\si_t^2}dY_t \ms\\
&=& [f_t-k_t]V_tdt +\sigma^1_{t}d B^1_{t}-l_tdY_t, \qq\q X_0-P_0=v-v_0.
\eeaa
It then follows that the identity (\ref{vmp}) holds $\hQ$-almost surely, and hence $\hP$-almost surely.

Similarly, applying
the constant variation formula for the linear SDE (\ref{p2}) and noting (\ref{sde4.1}) we obtain that, with $\phi_2(t,r)=
\exp(\int_r^t(f_u+g_u)du)$, for $0\le r,t\le T$,
\bea
\label{p3}
P_t&=&\phi_2(t,0)\Big\{v_0+\int_0^t\phi_2(0,r)h_rdr+\int_0^t\phi_2(0,r)\frac{H_rS_r}{(\si_r^2)^2}dY_r\Big\}\\
&=&\phi_2(t,0)\Big\{v_0+\int_0^t\phi_2(0,r)h_rdr+\int_0^t\phi_2(0,r)\frac{H_rS_r}{(\si_r^2)^2}[H_r(X_r-P_r)dr+\sigma_r^2dB_r^2]\Big\}. \nonumber
\eea
Now plugging (\ref{vmp}) into  (\ref{p3}), and applying Fubini, we have
\bea
\label{P}
P_t&=&\phi_2(t,0)\Big\{v_0+\int_0^t\phi_2(0,r)h_rdr+(v-v_0)\int_0^t\phi_1(r,0)\phi_2(0,r)k_rdr \nonumber \\
&&+\int_0^t\phi_1(0,r)\sigma_r^1\int_r^t\phi_1(u,0)\phi_2(0,u)k_u du dB_r^1\nonumber \\
&&+\int_0^t[\phi_2(0,r)l_r-\phi_1(0,r)l_r\int_r^t\phi_1(u,0)
\phi_2(0,u)k_udu]dB_r^2\Big\}\\
&=&\phi_2(t,0)\Big\{v_0+\int_0^t\phi_2(0,r)h_rdr+(v-v_0)\phi_3(t,0)\nonumber \\
 &&+\int_0^t\phi_1(0,r)\sigma_r^1\phi_3(t,r)dB_r^1+\int_0^t[\phi_2(0,r)l_r-\phi_1(0,r)\phi_3(t,r)l_r]dB_r^2\Big\}, \nonumber
 \eea
where $\phi_3(t,r)\dfnn\int_r^t\phi_1(u,0)\phi_2(0,u)k_udu$. This proves (\ref{p1}), whence the proposition.
\qed

\subsubsection{\bf  Unbounded Coefficients Case}

We note that Theorem \ref{well-posed} as well as the discussion so far rely heavily on the
assumption that all the coefficients are bounded, especially
$H$ and $G$ (see Assumption \ref{assump2}). However, in our applications we will see that  the coefficients
$H=-G=\beta$, where $\beta$ is the insider trading intensity which, at least in the optimal case, will 
satisfy $\lim_{t\to T^-}\beta_t=+\infty$, violating  Assumption \ref{assump2}.
In other words, the closed-loop system will exhibit a certain Brownian ``bridge" nature (see also, e.g., \cite{B98, CCD11,CCD13a}),
for which the well-posedness result of Theorem \ref{well-posed} actually does not apply.

To overcome such a conflict, we  introduce the following relaxed version of Assumption \ref{assump2}.
\begin{assum}
\label{assum20}
There exists a sequence $\{T_n\}_{n\ge1}$, with $0<T_n\nearrow T$, and a sequence of probability measures $\{\hQ^n\}_{n\ge 1}$ on
$(\Om^0, \cF^0)$, satisfying

\ms
(i) Assumption \ref{assump2} holds  for each $(\Om^0,\cF^0, \hQ^n)$ over $[0, T_n]$, $n\ge 1$

\ss
(ii) $\hQ^{n+1}\big|_{\cF^0_{T_n}} =\hQ^n$, $n\ge 1$.
\qed
\end{assum}

We shall refer to the sequence of probability measures ${\cQ}^0:=\{\hQ^n\}_{n\ge 1}$ as the {\it reference family of probability measures}, and the
associated sequence $\{T_n\}_{n\ge 1}$ as  the {\it announcing sequence}. Clearly, if the reference measure $\hQ^0$ exists, then $\hQ^n=\hQ^0|_{\cF^0_{T_n}}$, $n\ge 1$. It is known, however,  that in the dynamic observation case the Kyle-Back equilibrium may only exist on 
$[0,T)$ (see, e.g. \cite{CCD13b} and the references cited therein). In such a case the reference
family would play a fundamental role. A reasonable extension of the notion of $\hQ^0$-weak solution over $[0,T)$ is as follows.
\begin{defn}
\label{weaksol1}
Let ${\cQ}^0$ be a reference family of probability measures, with announcing sequence $\{T_n\}$. A sequence $\{(\Om^0, \cF^0, \hP^n, X^n, Y^n, B^{1, n}, B^{2, n})\}_{n\ge 1}$ is called a ${\cQ}^0$-weak solution of (\ref{sde4.1}) on $[0,T)$ if
for each $n\ge 1$, $(\Om^0, \cF^0, \hP^n, X^n, Y^n, B^{1,n}, B^{2,n})$ is a $\hQ^n$-weak solution on $[0,T_n]$.
\qed
\end{defn}

It is worth noting that if the coefficients of CMFSDE (\ref{sde4.1}) satisfy Assumption \ref{assump2} on each sub-interval $[0,T_n]$, then one
can apply Theorem \ref{well-posed} for each $n$ to get a ${\cQ}^0$-solution. Furthermore, since the solutions will be pathwisely
unique under each $\hQ^n$ over $[0,T_n]$, it is easy to check that $(X^{n+1}_t, Y^{n+1}_t, B^{1,n+1}_t, B^{2,n+1}_t)=(X^n_t, Y^n_t, B^{1,n}_t, B^{2,n}_t)$, $t\in [0, T_n]$, $\hQ^n$-a.s. We can then define a process $(X, Y, B^1, B^2)$ on $[0,T)$ by simply setting $(X_t, Y_t, B^{1}_t, B^{2}_t)=(X^n_t, Y^n_t, B^{1,n}_t, B^{2,n}_t)$, for $t\in [0, T_n]$, $n\ge 1$, and we shall refer to such a process as the {\it ${\cQ}^0$-solution} on $[0,T)$. The {\it ${\cQ}^0$-pathwise uniqueness on $[0,T)$} can be defined in an obvious way.
We have the following extension of Theorem \ref{well-posed}, whose proof is left for the interested reader.
\begin{thm}
\label{well-posed1}
Assume that Assumption \ref{assum20} is in force, and let ${\cQ}^0$ be the family of reference measures with announcing sequence $\{T_n\}$.
Assume further that Assumption \ref{assump2} holds for each $\hQ^n$ on $[0,T_n]$. Then CMFSDE (\ref{sde4.1}) possesses a ${\cQ}^0$-weak
solution on $[0,T)$, and it is ${\cQ}^0$-pathwisely unique on $[0,T)$.
\qed
\end{thm}

\section{Proof of Theorem \ref{well-posed}}
\setcounter{equation}{0}

In this section we prove Theorem \ref{well-posed}. We begin by making the following reduction: it suffices to consider the
SDE (\ref{sde4.1}) where  the initial state $X_0=v\equiv v_0$ is deterministic, that is, $s_0=0$.  Indeed, suppose that $(X^x, Y^x)$
is a weak solution
of (\ref{sde4.1}) along with some probability space $(\Om, \cF, \hP)$ and $\hP$-Brownian motion $(B^1, B^2)$, and $v$ is any  random variable defined on $(\hR, \sB(\hR))$, with normal distribution $\m\dfnn N(v_0, s_0)$,  we define the product space
$$\tilde\Om\dfnn \Om\otimes \hR, \q \tilde\cF\dfnn\cF\otimes\sB(\hR), \q\tilde \hP\dfnn \hP\otimes \m,
$$
and write generic element of $\tilde\o\in\tilde \Om$ as $\tilde\o=(\o, x)$. Then for each $t\ge 0$, the mapping $\tilde\o\mapsto X^x_t(\o)$ defines
a random variable on $(\tilde \Om, \tilde\cF, \tilde \hP)$, and $x\mapsto X_0^x\dfnn v(x)$ is a normal random variable with distribution $N(v_0,s_0)$ and is independent of $(B^1, B^2)$, by definition.  Bearing this in mind, throughout the section we shall assume
that the initial state $X_0=x$ is deterministic.

\subsection{Existence}

Our main idea to prove the existence of the weak solution is to ``decouple" the state and observation equations in
(\ref{sde4.1}) by considering the dynamics of the filtered state process $P_t \dfnn \hE^\hP [X_t|\cF^Y_t]$, $t\ge 0$, which is
known to satisfy an SDE, thanks to linear (Kalman-Bucy) filtering theory.

To be more precise, we consider the following system of SDEs on the reference probability space $(\Om, \cF, \hQ^0)$, on
which $(B^1, Y)$ is a Brownian motion:
\bea
\label{sdeQ0}
\left\{\ba{llll}
dX_t=[f_t X_t+g_t P_t +h_t]dt+ dB^1_t , \qq &X_0=x; \ms \\
dB^2_t= dY_t-[H_tX_t+G_t P_t]dt, &B^2_0=0; \ms\\
dP_t=[(f_t+g_t) P_t +h_t ]dt+S_t H_t \{dY_t -[H_t+G_t] P_tdt\}, & P_0 =x; \ms\\
dS_t=[2f_t S_t-H_t^2 S_t^2 +1]dt, & S_0= 0.
\ea\right.
\eea
We note that by Assumption \ref{assump2}, all coefficients $f, g, H, G$ are deterministic and $h\in L^2_{\hF^Y}(C([0,T]))$,
it is  easy to see that the linear system (\ref{sdeQ0}) has a (pathwisely) unique solution $(X_t, B^2_t, P_t)$ on $(\Om, \cF, \hQ^0)$.

Now let $L=\{L_t\}_{t\ge 0}$ be the solution to {the SDE}:
\be
\label{L}
d L_t = L_t(H_t X_t + G_t P_t ) dY_t, \q L_0 =1 .
\ee
Then $L$ is a  positive $\hQ^0$-local martingale, hence a
$\hQ^0$-supermartingale with $\hE^{\hQ^0}[L_t] \le L_0=1$. Furthermore, $L$ can be written as the Dol\'eans-Dade exponential:
\bea
\label{expL}
L_t=\exp\{\int_0^t (H_sX_s+G_sP_s)dY_s-\frac12 \int_0^t| H_sX_s+G_sP_s|^2ds\Big\},  \q t\in [0,T].
\eea

We have the following lemma.
\begin{lem}
\label{lemMG1}
Assume {that Assumption \ref{assump2} holds, and further that $h$ is bounded}. Then the process $L=\{L_t;t\ge0\}$ is a true $(\hF, \hQ^0)$-martingale on $[0,T]$.
\end{lem}
{\it Proof.} We follow the idea of that in \cite{Bens}.  Since $L$ is a supermartingale with $\hE^{\hQ^0}[L_t]\le 1$, we need only show that $\hE^{\hQ^0}[L_t]=1$, for
all $t\ge 0$.
To this end,  we define, for any $\e > 0$,
\beaa
\label{Le}
L_t^{\e}\dfnn\frac{L_t}{1+\e L_t}, \qq t\in[0, T].
\eeaa
Then clearly $0\le L^\e_t\le L_t\wedge \frac1\e$, and an easy application of It\^o's formula shows that
\be
\label{SDELe}
d L_t^{\e}=-\frac{\e L_t^2[H_t X_t +G_t P_t]^2}{(1+\e L_t)^3} dt +\frac{L_t[H_t X_t +G_t P_t]}{(1+\e L_t)^2} dY_t,\q t\ge0; \q L^\e_0=\frac1{1+\e}.
\ee
Since for each fixed $\e>0$,
\beaa
\Big|\frac{L_t[H_t X_t +G_t P_t]}{(1+\e L_t)^2} \Big|^2=\Big|\frac{\e L_t[H_t X_t +G_t P_t]}{\e(1+\e L_t)^2} \Big|^2\le\frac{[H_t X_t +G_t P_t]^2}{\e},
\eeaa
we see that the stochastic integral on the right hand side of (\ref{SDELe}) is a true martingale. It then follows that
\be
\label{ELe}
\hE^{\hQ^0}[L_t^{\e}]=\frac1{1+\e}-\hE^{\hQ^0}\Big[ \int_0^t \frac{\e L_t^2[H_t X_t +G_t P_t]^2}{(1+\e L_t)^3} dt\Big].
\ee
Next, we observe that $L_t >0$, and
\beaa
0 \leq \frac{\e L_t^2[H_t X_t +G_t P_t]^2}{(1+\e L_t)^3} = \frac{(\e L_t)L_t[H_t X_t +G_t P_t]^2}{(1+\e L_t)^3}\le L_t[H_t X_t +G_t P_t]^2.
\eeaa
Note that $L^\e$ is bounded. By sending $\e\to 0$ on both sides of (\ref{ELe}) and applying Dominated Convergence Theorem
we can then conclude that  $\hE^{\hQ^0}[L_t]=1$, provided
\bea
\label{bound}
\hE^{\hQ^0}\Big[ \int_0^T L_t[H_t X_t +G_t P_t]^2 dt\Big] < \infty.
\eea

It remains to check (\ref{bound}). To this end, let us define $X_t = X^1_t + \a_t$, where
\bea
\label{X1}
\left\{\ba{lll}
d\a_t=f_t\a_tdt+dB^1_t, \qq\qq & \a_0=x,  \ms \\
dX_t^1=[f_t X_t^1+g_t P_t +h_t]dt, & X_0^1=0.
\ea\right.
\eea
By Gronwall's inequality, it is readily seen that
\bea
\label{X1est}
|X^1_t| \leq C \int_0^t |g_s P_s +h_s|ds \leq C\Big[1+\int_0^t |P_s|ds\Big], \q t\in[0,T].
\eea
Here and in the sequel $C>0$ denotes a generic constant depending only on the bounds of the coefficients $f, g, H, G$, $h$, and
the duration $T>0$, which is allowed to vary from line to line.
Now, noting that $L_t$ is a super-martingale with $L_0=1$,
we deduce from (\ref{X1est}) that
\bea
\label{LX1est}
\hE^{\hQ^0}\Big[\int_0^T L_t|X^1_t|^2 dt\Big]& \leq& C\Big\{1+\hE^{\hQ^0}\Big[\int_0^T  \int_0^t L_t|P_s|^2 ds dt\Big]\Big\}\\
& \leq& C\Big\{1+\hE^{\hQ^0} \Big[ \int_0^T L_s|P_s|^2 ds\Big]\Big\}. \nonumber
\eea
Consequently we have
\bea
\label{bound1}
\hE^{\hQ^0}\Big[ \int_0^T L_t[H_t X_t +G_t P_t]^2  dt\Big] &\leq& C \hE^{\hQ^0}\Big[ \int_0^T L_t \big[|X^1_t|^2+|\a_t|^2+|P_t|^2 \big] dt\Big]\\
&\le&C\Big\{1+\hE^{\hQ^0}\Big[ \int_0^T L_t (|\a_t|^2+|P_t|^2 )dt\Big]\Big\}. \nonumber
\eea

Continuing, let us recall that the processes $P$ and $\a$ satisfy  (\ref{sdeQ0}) and (\ref{X1}), respectively. By It\^o's formula we see that
\bea
\label{Pa2}
\left\{\ba{lll}
d|\a_t|^2=[2f_t|\a_t|^2+1]dt+2\a_tdB^1_t;\ms\\
d|P_t|^2=\big[2M_t|P_t|^2+2P_th_t+S_t^2H_t^2 \big]dt +2S_tH_tP_tdY_t,
\ea\right.
\eea
where  $M_t \dfnn (f_t+g_t)- S_tH_t(H_t+G_t)$, $t\ge 0$.
Next, we define, for $\d>0$ and $t\in[0,T]$,
$$X_t^{\d}=\dfrac{X_t}{[1+\d |X_t|^2]^{1/2}}; \qq P_t^{\d}=\dfrac{P_t}{[1+\d |P_t|^2]^{1/2}}.
$$
Then $|X^\d_t|\le |X_t|\wedge \d^{-\frac12}$ and $|P^\d_t|\le |P_t|\wedge \d^{-\frac12}$, $\forall t$; and it is not hard to show
that  $\lim_{\d\to0}X^\d=X$, $\lim_{\d\to0}P^\d=P$, uniformly on $[0,T]$, in probability.
Now, define
\bea
\label{Ld}
dL_t^{\d}=L_t^{\d}[H_tX_t^{\d}+G_t P_t^{\d}]dY_t, \qq L_0^{\d}=1.
\eea
Since $X^\d$ and $P^\d$ are now bounded, $L^{\d}$ is a martingale and $\hE^{\hQ^0} [L_t^{\d}]=1$, $t\in [0, T]$.
Furthermore, by the stability of SDEs one shows that, possibly along a subsequence, $L^\d_t$ converges to $L_t$,
$\hQ^0$-a.s., $t\in[0,T]$.

Noting (\ref{Pa2}) and applying It\^o's formula we have, for $t\in[0,T]$
\bea
\label{Lda}
L_t^{\d}|\a_t|^2=x^2+\int_0^tL_s^{\d}[2f_s|\a_s|^2+1]ds+\int_0^t2L_s^{\d}\a_sdB_s^1+\int_0^tL_s^{\d}|\a_s|^2[H_sX_s^{\d}+G_sP_s^{\d}]dY_s. \nonumber
\eea
Since $\a$ has finite moments for all orders (see (\ref{X1})), the boundedness of $X^\d$ and $P^\d$ then renders the two stochastic
integrals on the right hand side above both true martingales. Thus, taking expectations on both sides above, and
applying Gronwall's inequality, we get
\bea
\label{Ldaest0}
\hE^{\hQ^0} \big[L_t^{\d}|\a_t|^2\big] \leq C, \qq \forall t\in[0,T],
\eea
where $C$ is a constant independent of $\d$. Applying Fatou's Lemma  we then get that
\bea
\label{La2est}
\hE^{\hQ^0}\big[ L_t|\a_t|^2\big] \leq \liminf_{\d\to0}\hE^{\hQ^0}[L^\d_t|\a_t|^2]\le C.
\eea
%
%
Finally, noting (\ref{Pa2}) and applying It\^o's formula again we have
\bea
\label{LdP2}
d L_t^{\d}|P_t|^2 &=&L_t^{\d}\big[2M_t|P_t|^2+2P_th_t +S_t^2H_t^2 \big]dt +2S_tH_t L_t^{\d}P_tdY_t\\
&& +L_t^{\d}|P_t|^2[H_tX_t^{\d}+G_t P_t^{\d}]dY_t+2S_tH_t L_t^{\d}P_t[H_t X_t^{\d}+G_t P_t^{\d}]dt. \nonumber
\eea
By similar arguments as before, and noting that $|X^\d_t|\le |X_t|$ and $|P^\d|\le |P_t|$, one shows
\beaa
\hE^{\hQ^0}[ L_t^{\d}|P_t|^2] &\leq& C\Big\{1+\hE^{\hQ^0}\Big[ \int_0^t L_s^{\d}[|P_s|^2+|X_s|^2]ds\Big]\Big\}\\
&\le & C\Big\{1+\hE^{\hQ^0}\Big[ \int_0^t L_s^{\d}[|P_s|^2+|X^1_s|^2+|\a_s|^2]ds\Big]\Big\}, \q t\in[0,T].
\eeaa
This, together with  (\ref{LX1est}), implies that
\beaa
\label{LdP2est}
\hE^{\hQ^0}[ L_t^{\d}|P_t|^2 ]\leq C\Big\{1+\hE^{\hQ^0}\Big[\int_0^t L_s^{\d}[|P_s|^2+|\a_s|^2]ds\Big]\Big\}, \q t\in[0,T].
\eeaa
Applying the Gronwall inequality and recalling (\ref{Ldaest0}) we then obtain
\bea
\label{SDE}
\hE^{\hQ^0}[ L_t^{\d}|P_t|^2] \leq C\Big\{1+\hE^{\hQ^0}\Big[ \int_0^t L_s^{\d}|\a_s|^2ds\Big]\Big\}\le C, \q t\in[0,T].
\eea
By Fatou's lemma one again shows that $\hE^{\hQ^0}[ L_t|P_t|^2] \leq C$, for all $t\in[0. T]$. This, together with (\ref{La2est})
and (\ref{bound1}), leads to (\ref{bound}). The proof is now complete.
\qed

We can now complete the proof of existence. Since $L_t$ is a $(\hF , \hQ^0)$ martingale, we  define a probability measure $\hP$ by $\frac{d\hP}{d\hQ^0}\big|_{\cF_T}=L_T$, and apply the Girsanov Theorem so that
$(B^1, B^2)$ is a $\hP$-Brownian motion on $[0,T]$. Now, by looking at the first two equations of
(\ref{sdeQ0}), we see that $(\Om, \cF, \hP, X, Y, B^1, B^2)$ would be a weak solution to (\ref{sde4.1}) if we can show that
\bea
\label{Peq}
P_t = \hE^{\hP} [X_t|\cF^Y_t], \q t\in[0, T], \q \hP\mbox{-a.s.}
\eea

To prove (\ref{Peq}) we proceed as follows. We  consider the following linear filtering problem on the space $(\Om, \cF, \hP)$:
\bea
\label{filter}
\left\{
\ba{lll}
d\a_t=f_t\a_tdt+dB^1_t, \qq &\a_0=x;\\
d\b_t=H_t \a_tdt +dB_t^2, &\b_0=0.
\ea
\right.
\eea
%
%
Denote $\widehat\a_t = \hE^\hP[\a_t|\cF^\b_t]$, $t\ge0$.  Then by linear filtering theory, we know that  $\widehat\a$ satisfies the
following SDE:
\be
\label{Kalman}
d\widehat\a_t=f_t \widehat\a_tdt+S_t H_t \{d\b_t - H_t \widehat\a_tdt\}, \qq \widehat\a_0=x,
\ee
where $S_t$ satisfies $(\ref{SDE4.21})$ (or (\ref{sdeQ0})). On the other hand, from (\ref{filter}) we see that $\a$ is $\hF^{B^1}$-adapted,
and  from (\ref{sdeQ0}) we see that $P$ is $\hF^Y$-adapted, therefore  we can apply Lemma \ref{lem4.1} to conclude
that   $\hE^\hP[\a_t|\cF_t^Y]=\hE^\hP[\a_t|\cF^\b_t]=\widehat\a_t$, $t \in[0,T]$.

Now let us define $\widetilde{P_t}=\hE^\hP [X_t|\cF^Y_t]$, $t\ge0$. Recall that $X=X^1+\a$, where $X^1$ satisfies a randomized
ODE (\ref{X1}) and  is obviously $\hF^Y$-adapted, we see that $\widetilde{P}=X^1+\widehat\a$,  and
it  satisfies the SDE:
\bea
\label{tildeP}
d\widetilde P_t=[f_t \widetilde P_t+g_t P_t +h_t]dt+S_tH_t\{d\b_t-H_t\widehat\a_t dt\}, \q t\ge0;  \qq \widetilde P_0=x.
\eea
%
%
%
Note that $X=X^1+\a$ and $\tilde P=X^1+\widehat \a$ we see that
\beaa
\label{innovation}
 d\b_t - H(t)\widehat\a_t dt=H_t (\a_t- \widehat\a_t )dt+dB^2_t= dY_t-[H_t\widetilde P_t+G_t P_t]dt.
\eeaa
Then (\ref{tildeP}) implies that
\be
\label{tildeP1}
d\widetilde P_t=[f_t\widetilde P_t+g_t P_t + h_t]dt+S_tH_t\{dY_t-[H_t\widetilde P_t+G_t P_t]dt\}.
\ee
Define $\D P_t=P_t-\widetilde P_t$, then it follows from (\ref{sdeQ0}) and (\ref{tildeP1}) that
\beaa
\label{SDE}
d\D P_t=[f_t -S_t H_t^2]\D P_tdt, \qq \D P_0=0.
\eeaa
Thus $\D P_t \equiv 0$, $\hP$-a.s., for any $t\ge 0$. That is, (\ref{Peq}) holds, proving the existence.

It is worth noting that the weak solution that we have constructed is actually a $\hQ^0$-weak solution.



\subsection{Uniqueness}

Again we need only consider the solutions with deterministic initial state. We first note that if $(\Om^0, \cF^0, \hP, \hF^0, X, Y, B^1, B^2)$
is a $\hQ^0$-weak solution to (\ref{sde4.1}), then we can assume without loss of generality  that $\hF^0=\hF^{B^1, B^2}$, hence Brownian.
Next, we can define $\widetilde{P_t}=\hE^\hP [X_t|\cF^Y_t]$. We are to show that $\widetilde{P_t}$ satisfies an SDE of the
form as that in (\ref{sdeQ0}) under $\hQ^0$, from which we shall derive the $\hQ^0$-pathwise uniqueness.

To this end, we recall that, as a $\hQ^0$-weak solution, one has $\hP\sim \hQ^0$.
Define a
$\hP$-martingale $Z_t\dfnn \hE^{\hP}\big[\frac{d\hQ^0}{d\hP}\big|\cF_t\big]$, $t\ge 0$. Since $(B^1,Y)$ is a $\hP$-semi-martingale
with decomposition:
\bea
\label{Ydec}
B^1_t=B^1_t; \q  Y_t=\int_0^t (H_sX_s+G_s\widetilde P_s)ds +B^2_t,\qq t\ge0.
\eea
By Girsanov-Meyer Theorem (see, e.g., \cite[Theorem III-20]{prot}), it is a $\hQ^0$-semi-martingale with the decomposition
$(B^1_t,Y_t)=(N^1_t,N^2_t)+(A^1_t, A^2_t)$, where $N=(N^1,N^2)$ is a $\hQ^0$-local martingale of the form
\beaa
 N^1_t= B^1_t-\int_0^t\frac1{Z_s}d[Z,B^1]_s,\q
 N^2_t = B^2_t-\int_0^t\frac1{Z_s}d[Z,B^2]_s, \q t\ge0,
\eeaa
and $A=(A^1,A^2)$ is a finite variation process. Since by assumption $(B^1,Y)$ is $\hQ^0$-Brownian motion, we have $A\equiv0$.
In other words, it must hold that
\bea
\label{B2Y}
B^1_t= B^1_t-\int_0^t\frac1{Z_s}d[Z,B^1]_s,\q Y_t=B^2_t-\int_0^t\frac1{Z_s}d[Z,B^2]_s, \q t\ge0.
\eea
Consider now a $(\hF, \hP)$-martingale $dM_t=Z^{-1}_tdZ_t$. Since $\hF$ is Brownian, applying Martingale Representation
Theorem we see that there exists a
process $\th=(\th^1, \th^2)\in L^2_{\hF}([0,T])$ such that $dM_t=\th^1_tdB^1_t+\th^2_t dB^2_t$, $t \in [0,T]$. Thus (\ref{B2Y})
amounts to saying that
$$ [M, B^1]_t=\int_0^t\th^1_sds\equiv 0;\q Y_t=B^2_t-[M,B^2]_t=B^2_t-\int_0^t\th^2_s ds, \q t\ge 0.
$$
Comparing this to (\ref{Ydec}) we have $\th^1\equiv0$ and  $\th^2\equiv -(HX+G\widetilde P)$. That is, $Z=L^{-(HX+G
\widetilde P)}$,
the solution to the SDE:
\beaa
dZ_t=Z_tdM_t=-Z_t (H_tX_t+G_t\widetilde P_t)dB^2_t, \q t\in[0,T], \q Z_0=1,
\eeaa
and hence it can be written as the Dol\'eans-Dade stochastic exponential:
\bea
\label{Zexp}
Z_t= \exp\Big\{-\int_0^t(H_sX_s+G_s\widetilde P_s)dB^2_s-\frac12\int_0^t |H_sX_s+G_s\widetilde P_s|^2ds\Big\}.
\eea

Let us now consider again the following filtering problem on probability space $(\Om, \cF, \hP)$.
\bea
\label{filter1}
\left\{
\ba{lll}
d\a_t=f_t\a_tdt+dB^1_t, \qq &\a_0=x.\\
d\b_t=[H_t\a_t ]dt +dB_t^2, &\b_0=0.
\ea
\right.
\eea
As before, we know that $\widehat\a_t = \hE^\hP[\a_t|\cF^\b_t]$ satisfies the SDE:
\be
\label{SDE10}
d\widehat\a_t=f_t\widehat\a_tdt+S_t H_t \{d\b_t - H_t\widehat\a_tdt\},  \qq \widehat\a_0=x,
\ee
where $S_t$ satisfies (\ref{SDE4.21}). Since  $L^{-(HX+G\widetilde P)}= Z$ is a $\hP$-martingale,
we can  apply  Lemma \ref{lem4.1} again to conclude that $\hE^\hP[\a_t|\cF_t^Y]=\hE^\hP[\a_t|\cF_t^\b]=\widehat\a_t$.

Now let $X_t=\a_t + X_t^1$, and $Y_t=\b_t+Y_t^1$ as before, where $X^1$ satisfies the ODE (\ref{X1}), and  $Y^1$ satisfies the ODE
\bea
\label{Y1}
dY_t^1=[H_tX_t^1+G_t \widetilde{P_t}] dt,\qq
Y_0^1=0.
\eea
Furthermore, since $X_t^1$ is $\hF^Y$-adapted, we have
$\widetilde{P_t}=X_t^1+\widehat\a_t$. Combining (\ref{X1}) and (\ref{SDE10}) we see that  $\widetilde{P_t}$ satisfies the  SDE:
\be
\label{SDE12}
d\widetilde P_t=[(f_t+g_t )\widetilde P_t +h_t]dt+S_tH_t\{d\b_t-H_t\widehat\a_t dt\}, \q \widetilde P_0=x.
\ee
Since $\b=Y-Y^1$, we derive from (\ref{Y1}) that
\beaa
 d\b_t - H_t\widehat\a_t dt= dY_t - dY_t^1- H_t\widehat\a_t dt= dY_t-[H_t\widetilde P_t+G_t \widetilde P_t]dt,
\eeaa
and (\ref{SDE12}) becomes
\be
d\widetilde P_t=[(f_t+g_t )\widetilde P_t +h_t]dt+S_tH_t\{dY_t-[(H_t+G_t )\widetilde P_t]dt\},\q t\in[0,T], \q \widetilde P_0=x.
\ee
That is, $\widetilde P_t$ satisfies the same SDE as $P_t$ does in (\ref{sdeQ0}) on the reference space $(\Om, \cF, \hQ^0)$.

To finish the argument, let  $(\Om, \cF, \hP^i, \hF, X^i, Y^i, B^{1,i}, B^{2,i})$, $i=1,2$ be any two $\hQ^0$-weak solutions,
and define $\widetilde P^i_t\dfnn \hE^{\hP^i} [X^i_t|\cF^{Y^i}_t]$, $t\ge0$, $i=1,2$. Then the arguments above show that
$(X^i, B^{2,i}, \widetilde P^i)$, $i=1,2$, are two solutions to the linear system of SDEs (\ref{sdeQ0}), under $\hQ^0$. Thus if
$(B^{1,1}, Y^1)\equiv (B^{1,2}, Y^2)$ under $\hQ^0$, then we must have $(X^1, B^{2,1}, \tilde P^1)\equiv
(X^2, B^{2,2}, \tilde P^2)$, under $\hQ^0$, which in turn shows, in light of  (\ref{Zexp}), that $\hP^1=\hP^2$. This proves
the $\hQ^0$-pathwise uniqueness of solutions to (\ref{sde4.1}).
\qed

\section{Necessary Conditions for Optimal Trading Strategy}
\setcounter{equation}{0}

In this section we study the optimization problem (\ref{j0}). We still denote the price dynamics observable by the insider
to be $V=\{V_t:t\ge0\}$, and assume that it satisfies the SDE:
\bea
\label{V1}
\left\{\ba{lll}
dV_t=[f_t V_t+g_t \hE^\hP[V_t|\cF^Y_t]+h_t]dt + \si^v_t dB^v_t, \q t\in[0,T], \ms\\
V_0=v\sim N(v_0,s_0);
\ea\right.
\eea
and we assume that the demand dynamics observable by the market makers, denoted by $Y=\{Y_t:t\ge0\}$, satisfies the
SDE
\bea
\label{Yeq1}
dY_t=[\b_t(V_t-\hE^\hP[V_t|\cF^Y_t])]dt +\si^z_t dB^z_t, \q t\in [0,T]; \qq Y_0=0.
\eea

We should note that in (\ref{V1}) and (\ref{Yeq1}) the probability $\hP$ should be understood as one defined on the
canonical space $(\Om^0, \cF^0, \hF^0)$, on which the solution to (\ref{V1}) and (\ref{Yeq1}) is $\hQ^0$-pathwisely
unique. For notational simplicity, from now on we shall denote $\hE=\hE^\hP$, when there is no danger of confusion.
 Moreover, note that
$$\hE[P_t(V_t-P_t)]=\hE[\hE[(V_t-P_t)P_t|{\cal F}_t^Y]]=0,$$
and that all the coefficients are now assumed to be deterministic, we can apply Proposition \ref{eu} to write
 the problem (\ref{j0}) as
\bea
\label{j1}
 J(\beta)&=&\int_0^T\beta_{t}\hE[(V_t-P_t)P_T]dt \nonumber\\
&= &\phi_2(T,0)\int_0^T\beta_t\phi_1(t,0)\Big{\{}s_0\phi_3(T,0)\\
 & & +\int_0^t[(l_r^2+(\sigma_r^v)^2)\phi^2_1(0,r)\phi_3(T,r)-\phi_1(0,r)\phi_2(0,r)l_r^2]dr\Big{\}}dt,  \nonumber
\eea
where $\phi_i$, $i=1,2,3$, and $l$, $k$ are defined by (\ref{phi123}) and (\ref{kl}), respectively, with $H=\b$. Now by integration by parts
we can easily check that
\beaa
&&\int_0^t[(l_r^2+(\sigma_r^v)^2)\phi^2_1(0,r)\phi_3(T,r)]dr=\int_0^t\phi_3(T,r)d[S_r\phi_1^2(0,r)]\\
&&\qq\qq=\phi_3(T,t)S_t\phi^2_1(0,t)-\phi_3(T,0)s_0+\int_0^t\phi_1(0,r)\phi_2(0,r)l_r^2dr,
\eeaa
we thus have
$J(\beta)=\phi_2(T,0)\int_0^T\beta_tS_t\phi_1(0,t)\phi_3(T,t)dt$,
and the original  optimal control problem (\ref{j0}) is equivalent to the following
\bea
\label{Jeq}
\sup_{\beta}\bar{J}(\beta) =\sup_{\b}\int_0^T\beta_tS_t\phi_1(0,t)\phi_3(T,t)dt.
\eea

Before we proceed any further let us specify the ``admissible strategy" and the standing assumptions on
the coefficients that will be used throughout this section. We note that the assumptions will be slightly stronger
than Assumption \ref{assump2}.
\begin{assum}
\label{assump3}
(i) All coefficients $f$, $g$, $h$, $\si^v$, and $\si^z$ are deterministic, continuous functions on $[0,T]$, such that
$\sigma_t^z \ge c$, $\sigma^v_t \ge c$ for all $t\in [0,T]$ for some constant $c>0$;

(ii) the trading intensity $\b$ is continuous on $[0,T)$, $\beta_t >0$ for all $t\in [0,T)$,  and $\lim_{t\to T^-} \beta_t
>0$ exists (it may be $+\infty$). Consequently, $\underline{\b}\dfnn\inf_{t\in[0,T]}\b_t>0$.
\end{assum}


\begin{rem}
\label{remark5}
{\rm (i) In practice it is not unusual to assume that $\lim_{t\to T^-}\b_t=\infty$, which amounts to saying that
the insider is desperately trying to maximize the advantage of the asymmetric information (cf. e.g., \cite{ABO12}). We shall actually prove
that this is the case for the optimal strategy, provided the Assumption \ref{assump3}-(ii) holds. In what follows
we say that a trading intensity $\b$ is {\it admissible} if it satisfies Assumption \ref{assump3}-(ii). By a slight  abuse of notations we still denote all admissible trading intensities by $\sU_{ad}$.

(ii) For $\b\in\sU_{ad}$, the well-posedness of CMFSDEs (\ref{V1}) and (\ref{Yeq1}) should be understood in the sense
of Theorem \ref{well-posed1}, and we shall consider its (unique) ${\cQ}^0$-solution.
\qed}
\end{rem}


We note that the solution of the Riccati equation (\ref{Riccati}) $S$, as well as the functions  $\phi_1$ and $\phi_3$
defined by (\ref{phi123}), depends on the choice of trading intensity $\b$. We shall at times denote them by $S^\b$,
 $\phi^\b_1$, and $\phi^\b_3$, respectively, to emphasize their dependence on $\b$. The following lemma is simple but useful  for our analysis.
\begin{lem}
\label{Sbeta}
Let Assumption \ref{assump3} be in force. Then  for any $\beta\in \sU_{ad}$,

(i) the  Riccati equation (\ref{Riccati})
has a solution $S=S^\b$ defined on $[0,T)$, such that $S^\b_t >0$, for all $t\in [0,T)$. Furthermore, there exists
a constant $C_\b>0$, depending on the bounds of the coefficients and $\underline{\b}$ in Assumption \ref{assump3}, such that $S^\b_t\le s_0+C_\b t$, $t\in [0,T]$;

\ss

(ii) $\|S^\b\|_\infty\le e^{KT}(s_0+KT)$, where $K\dfnn \|\si^v\|_\infty^2+2\|f\|_\infty$;

 \ss
 (iii) 
  $S^\b_{T^-}\dfnn\lim_{t\to T^-} S^\b_t<\infty$, that is, the solution $S^\b$ can be extended continuously
 to $[0, T]$;


(iv) $S^\b_ {T^-}=0$ if and only if $\lim_{t\rightarrow T^-}\phi_1^\b(t,0)=0$.
\end{lem}

{\it Proof.} Let $\b\in\sU_{ad}$ be given, and denote $S=S^\b$ and $\phi_1=\phi^\b_1$, etc., throughout the proof for simplicitly.

(i) First note that if $S_t =0$ for some $t<T$, we define $\t=\inf\{t\in[0,T); S_t =0\}$.
Then, from (\ref{Riccati}) we see that at $\t$ it holds that $\frac{dS_t}{dt}|_{t=\t}=|\sigma_{\t}^v|^2 >0$. But on the other hand by definition
of $\t$ we must have $\frac{S_{\t}-S_{\t-h}}{h} =-\frac{S_{\t-h}}{h}< 0$ for $h>0$ small enough, a contradiction.
That is,  $S_t >0$, $\forall t\in[0,T)$.

Next, let us denote, for $(t,s,\b)\in[0,T]\times (0,\infty)\times [0, \infty)$, the right side of (\ref{Riccati}) by
\beaa
G(t,s, \b)\dfnn (\si_t^v)^2+2f_t s-\Big[\frac{\beta s}{\si_t^z}\Big]^2.
\eeaa
Then  for any $\b\in\sU_{ad}$, it holds that
\bea
\label{Gb1}
G(t,s, \b_t)=-\frac{\beta^2_t}{|\sigma^z_t|^2}\Big[s-\frac{f_t|\sigma^z_t|^2}{\beta^2_t}\Big]^2+\frac{f^2_t|\sigma^z_t|^2}{\beta^2_t}+|\sigma^v_t|^2\leq \max_{t\in[0,T]}\Big\{ \frac{f^2_t|\sigma^z_t|^2}{{\underline\beta}^2}+|\sigma^v_t|^2 \Big\} \dfnn C_\b,
\eea
thanks to Assumption \ref{assump3}. Thus $S_t \leq s_0 + C_\b t$, for all $t\in[0,T]$, proving (i).

(ii) To find the bound that is independent of the choice of $\b$, we note that
$$ \frac{dS_t}{dt}=G(t, S_t, \b_t)\le  (\si^v_t)^2+2|f_t|S_t\le K(1+S_t), \q  \forall t\in[0,T],
$$
 where $K\dfnn \|\si^v\|_\infty^2+2\|f\|_\infty$. Thus
the result then follows from the Gronwall's inequality.

\ms
(iii) Since $G$ is quadratic in $s$, and $lim_{s\to\infty}G(t,s,\b)=-\infty$, it is easy to see from (\ref{Gb1}) that, for any
given $\b\in\sU_{ad}$,
\bea
\label{Gb2}
\max_{t,s} G^+(t,s,\b_t)=\max_{t,s}G(t,s,\b_t)\le C_\b.
\eea
On the other hand, we write
\beaa
S_t-s_0=\int_0^t G(r, S_r,\b_r)dr=\int_0^t G^+(r, S_r,\b_r)dr -\int_0^t G^-(r, S_r,\b_r)dr\dfnn I^+(t)-I^-(t),
\eeaa
where $I^\pm$ are defined in an obvious way. Since $I^+(\cd)$ and $I^-(\cd)$ are monotone increasing, both limits $I^+(T^-)$ and $I^-(T^-)$ exist, which may be
$+\infty$. But (\ref{Gb2}) implies that $I^+(T^-)<\infty$, and by (i),
$ I^-(t)=I^+(t)-S_t +s_0<I^+(t)+s_0$, for all $t\in [0,T)$, we conclude
that $I^-(T^-)<\infty$ as well. That is, $S_{T^-}\ge 0$ exists.

\ms
(iv) We rewrite the equation (\ref{Riccati}) as follows (recall the definition of $\phi_1$ (\ref{phi123})),
\bea
\label{phi1}
S_t&=&\exp(\log S_t)=s_0\exp\Big\{\int_0^t\frac{dS_t}{S_t}\Big\}=s_0\exp\Big\{\int_0^t(2f_t+\frac{(\sigma_t^v)^2}{S_t}-\frac{\beta_t^2S_t}{(\sigma_t^z)^2})dt\Big\}\\
&=&s_0\phi_1(t,0)\exp\Big\{\int_0^t(f_t+\frac{(\sigma_t^v)^2}{S_t})dt\Big\}.\nonumber
\eea
 Thus, the result follows easily from (iii). This completes the proof.
%
\qed

In the rest of the section we shall try to solve the optimization problem (\ref{Jeq}). We first note that by definition
 the quantity $S$ and hence $\phi_1$  and $\phi_3$ all depend on the choice of trading intensity function $\b$. Therefore
 (\ref{Jeq}) is essentially a problem of calculus of variation. We shall proceed by first looking for the first order necessary
 conditions, and then find the conditions that are sufficient for us to determine the optimal strategy.

 To begin with, let us denote, for any differentiable functional $F:\hC([0,T])\mapsto \hC([0,T])$ and any $\beta, \xi\in\hC([0,T])$,
 the directional derivative of $F$ at $\beta$ in the direction $\xi$ by
\bea
\label{DFxi}
\nabla_\xi F(\b)_t=\frac{d}{dy}F(\beta_t+y\xi_t)|_{y=0}.
\eea

We first give some useful directional derivatives that will be used frequently in the sequel.
Recall the solution $S$ to the
Riccati equation (\ref{Riccati}), and the functions $\phi_1$ and $\phi_3$, defined by (\ref{phi123}). Note that they are all
functionals of the trading intensity $\b\in\hC([0,T])$.
\begin{lem}
\label{diff}
Let $\xi=\{\xi_t\}$ be an arbitrary continuous function on $[0,T]$. Then the following identities hold, provided all the directional
derivatives exist:

(i) $\nabla_\xi\beta_t=\xi_t$;

\ms
(ii) $\dis\nabla_\xi S_t=-\phi^2_1(t,0)\int_0^t\nabla_\xi\beta_r \alpha_rS_r\phi_1^2(0,r) dr$; where $\alpha_t=\frac{2\beta_tS_t}{(\sigma_t^z)^2}$;

\ms

(iii) $\dis \nabla_\xi\phi_1(t,0)=-\phi_1(t,0)\int_0^t[\nabla_\xi\beta_r\alpha_r+\nabla_\xi S_r\rho_r]dr$; where $\rho_t=\frac{\beta_t^2}{(\sigma_t^z)^2}$;

\ms
(iv) $\dis \nabla_\xi\phi_1(0,t)=\phi_1(0,t)\int_0^t[\nabla_\xi\beta_r\alpha_r+\nabla_\xi S_r\rho_r]dr$;

\ms

(v) $\dis \nabla_\xi \phi_3(T,t)=\int_t^T\{ \nabla_\xi\phi_1(r,0)\phi_2(0,r)k_r+\phi_1(r,0)\phi_2(0,r)[\nabla_\xi\beta_r\alpha_r
+\nabla_\xi S_r\rho_r]\}dr$.
\end{lem}

 \noindent{\it Proof.}  (i)  is obvious. (iii)--(v) follows directly from chain rule. We only prove (ii).

To see this, recall (\ref{Riccati}). We have
$$
S_t(\beta+y\xi)=s_0+\int_0^t(\sigma^v_{r})^2dr+\int_0^t\Big\{2f_rS_r(\beta+y\xi)-\Big[\frac{S_r(\beta+y\xi)(\beta_r+y\xi_r)}{\sigma^z_{t}}\Big]^2\Big\}dr,
$$
and thus
\beaa
\nabla_\xi S_t&=&\int_0^t\frac{d}{dy}\Big\{2f_rS_r(\beta+y\xi)-\Big[\frac{S_r(\beta+y\xi)(\beta_r+y\xi_r)}{\sigma^z_{r}}\Big]^2\Big\}
\Big|_{y=0}dr\\
&=&\int_0^t\Big\{2f_r\nabla_\xi S_r-\frac{2S_r\beta_r}{\sigma^z_{r}}\Big[\frac{\beta_r\nabla_\xi S_r}{\sigma^z_{r}}+\frac{S_r\xi_r}{\sigma^z_{r}}\Big]\Big\}dr.
\eeaa
Denote $S^\nabla \dfnn \nabla_\xi S$, we see that it satisfies an ODE
$$
\frac{d}{dt}S^\nabla_t=2\Big[f_t-\frac{\beta_t^2S_t}{(\sigma^z_{t})^2}\Big]S^\nabla_t-\frac{2\b_tS_t^2\xi_t}{(\sigma^z_{t})^2}=2[f_t-k_t]S^\nabla_t-\frac{2\b_tS_t^2\xi_t}{(\sigma^z_{t})^2},
\qq S^\nabla_0=0.
$$
Solving it and noting that $\xi=\nabla_\xi\b$ and $\alpha_t=\frac{2\beta_tS_t}{(\sigma_t^z)^2}$ we obtain
$$
S^\nabla_t=-\int_0^t\frac{2\beta_rS_r^2\xi_r}{(\sigma^z_{r})^2}
\exp\Big\{\int_r^t2(f_u-k_u)du\Big\}dr=-\phi^2_1(t,0)\int_0^t\nabla_\xi\b_r \alpha_rS_r\phi_1^2(0,r) dr,$$
proving (ii), whence the lemma.
\qed

We are now ready to prove the following necessary conditions of optimal strategies for the original control problem (\ref{j0}).
\begin{thm}
\label{op1}
Assume that Assumption \ref{assump3} is in force. Suppose that $\beta\in\sU_{ad}$ is an optimal strategy of the problem (\ref{j0}),  then

(1) it holds that
$$\frac{\phi_1(t,0)(\sigma_t^z)^2}{2\b_tS_t}\phi_3(T,t)+\frac{1}{S_t}\Phi_t
=\int_t^T[\b_r\phi_1(r,0)\phi_3(T,r)+\frac{\b_r^2}{(\sigma_t^z)^2}\Phi_r]dr
$$
where
$$
\Phi_t=\phi_1^2(t,0)\int_0^t\b_rs_r\phi_1(0,r)dr\int_t^Tg_r\phi_1(r,0)\phi_2(0,r)dr.
$$

(2) Furthermore, $\lim_{t\rightarrow T^-}\b_t= \infty$,    and  consequently,
$\lim_{t\rightarrow T^-}S_t= 0$. In particular,
$$
P_T=\hE^\hP[V_T|{\cal F}_T^Y]=V_T, \qq \hP\mbox{-a.s.}
$$
\end{thm}

\noindent{\it Proof.} Suppose that $\beta\in\sU_{ad}$  is an optimal strategy of the problem (\ref{j0}). Then it is also an optimal of the problem (\ref{Jeq}). Thus for any  function $\xi\in \hC([0,T])$, it holds that
$$
\nabla_\xi\bar{J}(\beta)=\frac{d\bar{J}(\b+y\xi)}{dy}\Big|_{y=0}=0,
$$
or equivalently,
\bea
\label{nec1}
0&=&\int_0^T[\nabla_\xi\beta_t S_t\phi_1(0,t)\phi_3(T,t)+\beta_t\nabla_\xi S_t\phi_1(0,t)\phi_3(T,t)+\beta_t
S_t \nabla_\xi\phi_1(0,t)\phi_3(T,t)\nonumber\ms\\
&&\qq+\beta_t S_t\phi_1(0,t)\nabla_\xi\phi_3(T,t)]dt.
\eea
Then substituting $\nabla_\xi \b, \nabla_\xi S$, and $\nabla_\xi\phi_i$, $i=1,2,3$ in
Lemma \ref{diff} into (\ref{nec1}) and changing the order of integration if necessary
we obtain that, formally,
\bea
\label{neop0}
\frac{\phi_1(t,0)(\sigma_t^z)^2}{2\b_tS_t}\phi_3(T,t)+\frac{1}{S_t}\Psi_t
=\int_t^T[\b_r\phi_1(r,0)\phi_3(T,r)+\frac{\b_r^2}{(\sigma_t^z)^2}\Psi_r]dr
\eea
where
\bea
\label{Psi}
\Psi_t=\phi_1^2(t,0)[\phi_1(t,0)\phi_2(0,t)-\phi_3(T,t)]\int_0^t\b_rS_r\phi_1(0,r)dr.
\eea
To justify the identity (\ref{neop0}) we now show that both sides of (\ref{neop0}) are finite for any $\b\in \sU_{ad}$ (note
that it is possible that $\beta_t\to \infty$, as $t\to T^-$). To this end, we first note that
$\phi_1(t,0)$ and $\phi_2(r,s)$ are bounded, and
\bea
\label{phi3}
\phi_3(T,t)=\phi_1(t,0)\phi_2(0,t)-\phi_1(T,0)\phi_2(0,T)-\int_t^Tg_r\phi_1(r,0)\phi_2(0,r)dr,
\eea
thus it is  also bounded, and clearly $\lim_{t\to T}\phi_3(T,t)=0$. Furthermore,
we rewrite
(\ref{Psi}) as
\bea
\label{Psi1}
\Psi_t=\phi_1^2(t,0)G(t, T)\int_0^t\b_rS_r\phi_1(0,r)dr,
\eea
where
$G(t,T)\dfnn \phi_1(t,0)\phi_2(0,t)-\phi_3(T,t)=\phi_1(T,0)\phi_2(0,T)+\int_t^Tg_r\phi_1(r,0)\phi_2(0,r)dr$,
thanks to (\ref{phi3}). We claim that the following integral
\bea
\label{RHS}
\int_0^T[\b_r\phi_1(r,0)\phi_3(T,r)+\frac{\b_r^2}{(\sigma_t^z)^2}\Psi_r]dr\dfnn I_1+I_2,
\eea
where $I_1$ and $I_2$ are defined in an obvious way, is well-defined.
Indeed,  Assumption \ref{assump3} and the boundedness of $S$ imply that, modulo a universal constant,
\bea
\label{I1}
\b_t\phi_1(t,0)\phi_3(T, t)\sim \b_t \exp\Big\{\int_0^t\Big(f_u-\frac{\beta_u^2 S_u}{|\sigma^z_u|^2}\Big)du\Big\}\sim \frac{\b_t}{\exp(\int_0^t \b_u^2du)}.
\eea
On the other hand, by (\ref{Psi1}) it is easy to see that
$
\b_t^2\Psi_t
\sim\b_t^2\phi_1^2(t,0)\int_0^t\b_rS_r\phi_1(0,r)dr$, and
\bea
\label{intI1}
\int_0^t\b_rS_r\phi_1(0,r)dr&=& \int_0^t \b_rS_r \exp\Big\{\int_0^r\Big[-f_u+\frac{\beta_u^2 S_u}{|\sigma^z_u|^2}\Big]du\Big\} dr\nonumber\\
&\leq& C \int_0^t \b_rS_r \exp\Big\{\int_0^r\frac{\beta_u^2 S_u}{|\sigma^z_u|^2}du\Big\} dr \\
&=& C \int_0^t \frac{|\sigma^z_r|^2}{\b_r} d\Big[\exp\Big\{\int_0^r\frac{\beta_u^2 S_u}{|\sigma^z_u|^2}du\Big\}\Big] \le
C\exp\Big\{\int_0^t\frac{\beta_u^2 S_u}{|\sigma^z_u|^2 }du\Big\}. \nonumber
\eea
Here in the above $C>0$ is a generic constant, which depends only on the bounds of the coefficients and
$\underline\b$ in Assumption \ref{assump3}, and is allowed to vary from line to line. Thus, similar to (\ref{I1}),
we derive from (\ref{Psi1}) and (\ref{intI1}) that
\bea
\label{bPsi}
\b_t^2\Psi_t \le C \b^2_t \phi_1(t,0)G(t, T)
\exp\Big\{\int_0^t\frac{\beta_u^2 S_u}{|\sigma^z_u|^2} du\Big\}\sim   \b^2_t \phi_1(t,0)\sim
\frac{\b^2_t}{\exp(\int_0^t \b_r^2dr)} .
\eea
But notice that $\b_t\le 1+\b^2_t$, and that for any $\d>0$, we have
\beaa
\int_0^{T-\d}\frac{\b^2_t}{\exp\{\int_0^t \b^2_udu\}}dt\le 1-e^{-\int_0^{T-\d}\b^2_udu}\le 1,
\eeaa
we can easily derive from (\ref{I1}) and (\ref{bPsi})  that both
$I_1$ and $I_2$ in (\ref{RHS}) are finite, proving the claim.

We can now use the identity (\ref{neop0})  to prove both conclusions of the theorem.
 We begin by observing that $\lim_{t\to T^-}S_t=0$ must hold. In fact,
multiplying $S_t$ on both sides of (\ref{neop0}) and then taking limits
$t\to T^-$, and noting that $\lim_{t\to T^-}\phi_3(T, t)=0$,  we conclude that
$\lim_{t\to T^-}\Psi_t=0$. But then
the equation (\ref{Psi}), together with the fact that $\lim_{t\to T^-}\phi_3(T,t)=0$ but
$\lim_{t\to T^-}\phi_2(0,t)\neq 0$, implies that
$\lim_{t\to T^-}\phi_1(t, 0)=0$, and hence $\lim_{t\rightarrow T^-}S_t= 0$, thanks to Lemma \ref{Sbeta}.

We now claim  that $\lim_{t\rightarrow T^-}\b_t= \infty$. Indeed, suppose not, then we have
$$ \lim_{t\to T^-}\frac{dS_t}{dt}=\lim_{t\to T^-} G(t, S_t, \b_t)=\lim_{t\to T^-}(\si^v_t)^2\ge c>0,
$$
which  is a contradiction, since $S_t>0$, $\forall t\in[0,T)$, and $\lim_{t\to T^-}S_t=0$, proving the claim.

Now note that $S_t$ is the variance of the process $V_t-P_t$, $t\ge 0$, the facts that $\lim_{t\to T^-} S_t=0$
and that both processes $V$ and $P$ are continuous lead to that $V_T=P_T$, $\hP$-a.s.. This completes
the proof of part (2).

%

It remains to prove part (1). But note that $\phi_1(T, 0)=0$, we see from (\ref{Psi1}) that $\Psi_t=\Phi_t$, as
$G(t, T)=\int_t^T g_r\phi_1(r, 0)\phi_2(0,r)dr$. Thus (1) follows directly from (\ref{neop0}).
\qed

%
\ms
\begin{rem}
{\rm   (i) Theorem \ref{op1} amounts to saying that, similar to the static information
case (cf. e.g., \cite{ABO12, Dani}), the optimal strategy $\beta$ for the optimization problem (\ref{j0}) should also maximize the
advantage on the asymmetry of  information near the terminal time $T$,  i.e., $\lim_{t\rightarrow T^-}\b_t= \infty$.
Also, in equilibrium, such asymmetry of the information will disappear at the terminal time since
$S_T=\lim_{t\rightarrow T^-}S_t=0$, that is,  $P_T=\hE^\hP[V_T|{\cal F}_T^Y]=V_T$, $\hP$-a.s.,
despite the fact that the insider only observes $V_t$ at time $t<T$ (see also, e.g., \cite{B98, CCD11, CCD13a, CCD13b}).

(ii) Although Theorem 5.5 gives only the necessary condition of the optimal strategy, the well-posedness of the optimal closed-loop system 
is guaranteed by Theorem \ref{well-posed1}. In other words, combining Theorem \ref{op1} and Theorem \ref{well-posed1} we have proved the existence of the Kyle-Back equilibrium for the problem (\ref{V})---(\ref{j0}) under our assumptions.
\qed }
\end{rem}

\section{Worked-out Cases and Examples}
\setcounter{equation}{0}

In general, it is not easy to find out an closed-form optimal strategy for the original problem (\ref{j0}), although we somehow
predicted its behavior in Theorem \ref{op1}. In this subsection we consider a special case for which the optimal strategy
can be found explicitly. We show that it does possess the properties that we presented in the last sections, which in
a sense justifies our results. More precisely, we shall consider a case where the market price does not impact
the underlying asset price, namely, $g_t\equiv 0$ in equation (\ref{v1}). Recall from Lemma \ref{diff}-(ii) that  $\alpha_t=\frac{2\b_tS_t}{(\sigma_t^z)^2}$, $t\in[0,T]$.

\begin{thm}
\label{geq0}
 Assume that all the assumptions of Theorem \ref{op1} are in force and further that
 the coefficient  $g_t\equiv 0$ in the dynamics (\ref{v1}). Define $\a_0\ge 0$ so that  
 $$\frac{\a_0^2}{4}\dfnn \frac{s_0+\int_0^T\phi^2_2(0,r)(\sigma^v_{r})^2dr}{\int_0^T(\sigma_t^z)^2dt}.$$
Then, the solution to the optimization
 problem (\ref{j0})  is given as follows:

 (i) the optimal strategy is given by $\dis \beta_t=\frac{\a_0\phi_2(t, 0)(\sigma^z_{t})^2}{2S_t}$;

 \ms
(ii) the error variance  of  price $P_t$
\bea
\label{optS}
 S_t=\phi^2_2(t,0)\int_t^T\Big(\frac{\a_0^2}{4}(\sigma^z_{r})^2-\phi^2_2(0,r)(\sigma^v_{r})^2\Big)dr, \q t\in[0,T];
 \eea

  (iii) the corresponding expected payoff is given by
  $$\dis J(\beta)=\frac{\a_0\phi_2(T,0)}{2}\int_0^T{(\sigma_t^z)^2}dt;
$$

 (iv) the market price is given by
 $$P_t=\hE^\hP[V_t|{\cal F}_t^Y]=\phi_2(t,0)\Big[v_0+\int_0^t\phi_2(0,r)h_rdr+\frac{\a_0}{2}Y_t\Big], \q t\in[0,T];
$$

\ms
(v) finally, it holds that $\dis \lim_{t\rightarrow T^-}\b_t=\infty$, $\lim_{t\rightarrow T^-}S_t=0$,
and in particular,
\bea
\label{terminal}
P_T=\hE^\hP[V_T|{\cal F}_T^Y]=V_T, \qq \hP\mbox{-a.s.}
\eea
\end{thm}

{\it Proof.}  Let $\b\in\sU_{ad}$ be an optimal strategy. Then by Theorem \ref{op1},  we should
have $ \lim_{t\rightarrow T^-}S_t=0$, and hence $\phi_1(T,0)=\lim_{t\rightarrow T^-}\phi_1(t,0)=0$,
thanks to Lemma \ref{Sbeta}.

Since  $g_t\equiv 0$, we have
\bea
\label{phi31}
\phi_3(T,t)=\phi_1(t,0)\phi_2(0,t)-\phi_1(T,0)\phi_2(0,T)=\phi_1(t,0)\phi_2(0,t);
\eea
and (\ref{neop0}) now reads
\bea
\label{neop1}
\frac{\phi^2_1(t,0)\phi_2(0,t)(\sigma_t^z)^2}{2\b_tS_t}=\int_t^T\b_r\phi^2_1(r,0)\phi_2(0,r)dr.
\eea
Now recall from Lemma \ref{diff}-(ii) that  $\alpha_t=\frac{2\b_tS_t}{(\sigma_t^z)^2}$,
we can rewrite (\ref{neop1}) as
\bea
\label{neop2}
\phi_{12}(t)={\a_t}\int_t^T\b_r\phi_{12}(r)dr, \qq t\in [0,T],
\eea
where $\phi_{12}(t)\dfnn \phi^2_1(t, 0)\phi_2(0,t)$.
Now differentiating  with respect to $t$ on both sides of (\ref{neop2}) we obtain that
\bea
\label{dphi12}
\frac{d\phi_{12}(t)}{dt}= -\a_t\b_t\phi_{12}(t)+\frac{d\a_t}{dt}\int_t^T\b_r\phi_{12}(r)dr=
 \Big[-\a_t\b_t+\frac1{\a_t}\frac{d\a_t}{dt}\Big]\phi_{12}(t).
\eea
Now since $g\equiv 0$, we see from (\ref{phi123}) that $\phi_1(t,0)=\exp(\int_0^t(f_r-k_r)dr)$, $\phi_2(0,t)=\exp(-\int_0^tf_rdr)$, where, by (\ref{kl}), $k_t=\frac{\beta_t^2S_t}{(\sigma_t^z)^2}$. We can easily compute that
$$ \frac{d\phi_{12}(t)}{dt}=
2[f_t-k_t]\phi_{12}(t)-f_t\phi_{12}(t)=[f_t-2k_t]\phi_{12}(t).
$$
Plugging this into (\ref{dphi12}) and noting that, by definition, $2k_t=\a_t\b_t$, we obtain that
\bea
\label{dalpha}
\frac{d\a_t}{dt}=f_t\a_t, \qq \a_0=\frac{2\b_0s_0}{(\si^z_0)^2}.
\eea
Solving the above ODE we get
\bea
\label{alpha}
 \frac{2\b_tS_t}{(\sigma_t^z)^2}=\a_t=\a_0\exp\Big\{\int_0^t f_rdr\Big\}=\a_0\phi_2(t,0).
\eea
Consequently we obtain that $\dis \b_t=\frac{\a_0\phi_2(t,0)(\sigma_t^z)^2}{2S_t}$, proving (i).

\ms
(ii) Note that $S$ satisfies the ODE (\ref{Riccati}) and
$\phi_2(0,t)=\phi^{-1}_2(t,0)=\a_0/\a_t$, where $\a$ satisfies  (\ref{dalpha}),  it is easy to check that
\beaa
\frac{d}{dt}(\phi_2^2(0,t)S_t)
&=&\phi^2_2(0,t)\Big[(\si^v_t)^2+2f_tS_t-\Big(\frac{\b_tS_t}{\si^z_t}\Big)^2\Big]+S_t\a_0^2 \frac{d}{dt}\Big[\frac1{\a^2_t}\Big]\\
&=&\phi^2_2(0,t)\Big[(\si^v_t)^2+2f_tS_t-\Big(\frac{\b_tS_t}{\si^z_t}\Big)^2\Big]-2f_tS_t \phi^2_2(0,t)\\
&=&\phi^2_2(0,t)(\sigma^v_{t})^2-\frac{\a_0^2(\sigma^z_{t})^2}{4}.
\eeaa

Integrating both sides from $t$ to $T$, and noting that $S_{T^-}=0$ and $\phi_2(0,t)=\phi_2^{-1}(t, 0)$, we derive
(\ref{optS}). Furthermore,
setting $t=0$ in (\ref{optS}), one can then solve
$$\frac{\a_0^2}{4}\dfnn \frac{S_0+\int_0^T\phi^2(0,r)(\sigma^v_{r})^2dr}{\int_0^T(\sigma_t^z)^2dt},
$$
proving (ii).

\ms

%

(iii) Combining the expression of $\b$, $S$, as well as equations (\ref{j1}) and (\ref{phi31}), the corresponding expected payoff is
\bea\nonumber
J(\beta)&=&\phi_2(T,0)\bar{J}(\b)=\phi_2(T,0)\int_0^T\beta_tS_t\phi_1(0,t)\phi_3(T,t)dt\\\nonumber
&=&\phi_2(T,0)\int_0^T\beta_tS_t\phi_1(0,t)\big[\phi_1(t,0)\phi_2(0,t)\big]dt\\\nonumber
&=&\phi_2(T,0)\int_0^T\big[\frac{\a_0\phi_2(t,0)(\sigma_t^z)^2}{2S_t}\big]S_t\phi_2(0,t)dt=\frac{\a_0\phi_2(T,0)}{2}\int_0^T{(\sigma_t^z)^2}dt.
\eea
%

\ms
(iv) It follows from (\ref{p3}) and (\ref{alpha}) that the market price follows the dynamics
$$
P_t=\hE^\hP[V_t|{\cal F}_t^Y]=\phi_2(t,0)\Big[v_0+\int_0^t\phi_2(0,r)h_rdr+\frac{\a_0}{2}Y_t\Big], \q t\in[0,T].
$$

\ms

(v) Finally, again by Lemma \ref{Sbeta}, one has $\lim_{t\rightarrow T^-}\b_t=\infty$ and $\lim_{t\rightarrow T^-}S_t=0$.
In particular, $P_T=\hE[V_t|{\cal F}_T^Y]=V_T$,  $\hP$-a.s.
\qed

We note that Theorem \ref{geq0} contains several previously known results as special cases. We list them as follows.

\ms
{\bf 1. Static case.}  In this case $V_t \equiv v$, where $v\sim N(v_0,S_0)$.
%


Setting $f\equiv0, g\equiv0, h\equiv0, \sigma^v \equiv0$ in Theorem \ref{geq0} we have
$$\phi_2(t,0)\equiv 1, \qq\mbox{ and} \qq \a_0 = \Big(\frac{4 S_0}{\int_0^T(\sigma_t^z)^2dt}\Big)^{1/2}=2\l,
$$
where  $\l \dfnn \sqrt{S_0}\{\int_0^T(\sigma_t^z)^2dt\}^{-1/2}$ is the so-called {\it price sensitivity} or {\it Kyle's $\l$} (cf. \cite{ABO12}).
The optimal strategy is given by
$$ \dis \beta_t=\frac{\a_0(\sigma^z_{t})^2}{2S_t} = \frac{\big(\int_0^T (\sigma^z_t)^2 dt \big)^{1/2} (\sigma^z_{t})^2}{S_0^{1/2}  \int_t^T(\sigma^z_{r})^2dr};$$
and the corresponding expected payoff is given by
$$ J(\beta)=S_0^{1/2} \Big(\int_0^T{(\sigma_t^z)^2}dt \Big)^{1/2}. $$

Furthermore, one can easily check that, for  $t\in[0,T]$, the corresponding market price is given by
$$P_t=\hE^\hP[V_t|{\cal F}_t^Y] =v_0+\l Y_t, $$  and the corresponding mean square error is
$$S_t=\frac{S_0 \int_t^T(\sigma^z_{r})^2dr}{\int_0^T (\sigma^z_t)^2 dt}.$$
In particular, $S_T =0$, which implies that  $ V_0 =P_T$, $\hP$-a.s. These results coincide with those of \cite{ABO12}.

If we further assume that
%
%
$T=1$ and $\sigma^z_t \equiv \sigma$, where $\sigma>0$ is a constant, then the
optimal trading intensity becomes
$$
\b_t = \frac{\sigma}{\sqrt{S_0} (1-t)}, \qq 0 \leq t <1,
$$
the corresponding expected payoff is $J(\b)=\sigma \sqrt{S_0}$, and the corresponding market price is
$ P_t = v_0 + \frac{\sqrt{S_0}}{\sigma} Y_t$,  $0 \leq t <T$. We recover the results of Back \cite{B92}.

\ms
{\bf 2. Long-lived information case.} We now compare our results with that of Back-Pedersen \cite{B98} (see also Danilova \cite{Dani}), in which the
 insider continuously observes the dynamics of $V_t$ that is assumed to be a martingale.

By setting $T=1$, $f\equiv0$, $g\equiv0$, $h\equiv0$ and $\sigma^z\equiv1$, Theorem \ref{geq0} implies $\frac{\a_0^2}{4}=1$,
assuming 
$S_0 \dfnn 1- \int_0^1 (\sigma^v_s)^2 ds$, and the optimal trading intensity is
$$
\b_t = \dis \frac{1}{1-t-\int_t^1 (\sigma^v_s)^2 ds}, \qq 0 \leq t <1.
$$
The corresponding expected payoff is $J(\b)=1$, and the corresponding market price is  $ P_t = v_0 + Y_t$,  $0 \leq t <T$.

\bs
\no{\bf Acknowledgment.} We would like to express our sincere gratitude to the anonymous referees for their careful reading of the manuscripts and very insightful suggestions, which helped us to improve the quality of the paper significantly.

\end{document}